# REVISITING FAZLY–WEI–XU POINTWISE ESTIMATES FOR THE FOURTH-ORDER HÉNON EQUATION VIA BERNSTEIN'S TECHNIQUE

QUỐC ANH NGÔ AND TRUNG NGUYEN

ABSTRACT. In a remarkable work published in *Analysis and PDE* **8** (2015) 1541–1563, M. Fazly, J. Wei, and X. Xu establish the following pointwise estimate

$$-\Delta u \geqslant \frac{2}{n-4}\frac{|\nabla u|^2}{u} + \sqrt{\frac{2}{p+1-\frac{8}{n(n-4)}}}|x|^{\frac{\sigma}{2}}u^{\frac{p+1}{2}} \quad \text{in } \mathbf{R}^n$$

for any bounded, positive, $C^4$-solution $u$ to the the fourth-order Hénon equation

$$\Delta^2 u = |x|^\sigma u^p \quad \text{in } \mathbf{R}^n$$

with $n \geqslant 5$, $\sigma \geqslant 0$ but implicitly near 0, and $p > (n+4+2\sigma)/(n-4)$. Their argument relies on a sophisticated iteration argument in the fashion of the standard Moser proof. In this paper, we provide an alternative approach which relies on the Bernstein technique via the maximum principle. The heart of our argument is a suitable choice of an auxiliary function allowing us not only achieving the same pointwise estimate, which still holds even for $p = (n+4+2\sigma)/(n-4)$, but also relaxing the boundedness of solutions. In addition, the pointwise inequality actually holds for any $\sigma \geqslant 0$. This estimate leads to an interesting consequence of the existence of conformal metric of $\mathbf{R}^n$ having positive scalar and $Q$ curvatures. We also show that the coefficient $2/(n-4)$ is sharp, hence providing the first answer to the question in the op. cit. paper. Our proof appears to be simpler and constructive.

## 1. INTRODUCTION

The goal of this paper is to present another application of the Bernstein technique to obtain local and global pointwise estimates for solutions to the fourth-order elliptic equation in the whole Euclidean space $\mathbf{R}^n$ with $n \geqslant 5$, given by

$$\Delta^2 u = |x|^\sigma u^p \quad \text{in } \mathbf{R}^n \tag{1.1}$$

with $p > 1$ and $\sigma \geqslant 0$. In the existing literature, (1.1) is often called the fourth-order Hénon equation, which is an analogy of the classical Hénon equation proposed in [**Hen73**], which has its root in astrophysics.

Born more than a hundred years ago, the Bernstein technique, introduced in the local case by Bernstein himself in his pioneering papers [**Ber06**, **Ber10**], remains a powerful tool for establishing derivative estimates, either in pointwise or in integral sense, for solutions of elliptic equations. This technique is done through the use of auxiliary and cut-off functions together with the maximum principle. Consequently, higher order derivatives of solutions can be estimated in terms of lower order ones, but very often the size of the reference domain should appear in these estimates as a price we pay. For the simplest example demonstrating this technique, we refer the reader to [CDV22, page 1598].

The Bernstein technique has now been studied and extended to many other situations. In addition, it has been used to tackle many questions. As a matter of fact, it is impossible to provide a comprehensive account of all the results related to the Bernstein technique. For the interested readers, we refer to [LU64, Oba71, Ser71, Lio80, GS81, Bar91, BV91, CDV22, GL25], to name a few, and the references therein for several typical examples and their subsequent influence often in PDE theory and geometric analysis.

While the Bernstein technique is quite general, which is capable of handling a wider class of equations, some specific techniques built upon it can be used to obtain better estimates because it ties to the structure of the equation being studied. A particular example traces back to a remarkable paper by Modica in which he proved that any bounded $C^3$-solution $u$ to the second-order equation

$$\Delta u = F'(u) \quad \text{in } \mathbf{R}^n$$

with the nonlinearity $F \geqslant 0$ satisfies the following global gradient bound

$$|\nabla u|^2 \leqslant 2F(u) \quad \text{in } \mathbf{R}^n. \tag{1.2}$$

See [Mod85]. To achieve the goal, Modica used the following auxiliary defined by

$$P(x) = |\nabla u(x)|^2 - 2F(u(x)),$$

also called a $P$-function, and showed that it satisfies $\mathcal{L}(P) \geqslant 0$ for the second-order elliptic operator $\mathcal{L}$ defined by

$$\mathcal{L}(P) = \Delta P - 2F'(u)\frac{\nabla u}{|\nabla u|^2} \cdot \nabla P.$$

From this and with help of cut-off functions the conclusion $P \leqslant 0$ follows from the maximum principle. To obtain (1.2), the boundedness assumption is essential because this implies the smallness of $|\nabla u|^2$, which eventually controls the sign of the function $P$. However, a priori gradient bounds for unbounded solutions are also available in the literature, at least for certain class of equations. For e.g., by the classical Bernstein technique applied for the case $\mathcal{L} = \Delta$, it was proved in [MW08] that any positive $C^2$-solution $u$ to

$$\Delta u = u^p \tag{1.3}$$

in $B_R$ with $p < 0$ enjoys the following local gradient bound

$$|\nabla u|^2 \leqslant Cu(x)^2 + u(x)^{p+1} \quad \text{in } B_R(0) \tag{1.4}$$

for some absolute constant $C > 0$ depending only on $R$. The estimate (1.4) can also be called a logarithmic gradient estimate because it can be rewritten in a simpler way as follows

$$|\nabla \log u|^2 \leqslant C + u(x)^{p-1} \quad \text{in } B_R(0).$$

Clearly, (1.4) is similar to (1.2), but the right hand side of (1.4) is not a primitive of the right hand side of (1.3) due to the presence of the quadratic term $u^2$. (It is now known that radial solutions to (1.3) are unbounded, see [NNP26] and the references therein.) It is worth noting that the local gradient estimate (1.4) for (1.3) can be extended to other elliptic equation such as

$$-\Delta u = u^p \quad \text{in } \mathbf{R}^n, \tag{1.5}$$

see [MW08, page 1064]; see also [Li91]. Equation (1.5) is often called the Lane–Emden equation, that has been the center object of an enormous number of studies in this field.

It is by now well-known that an estimate such as (1.2) implies Liouville type results and monotonicity properties of the relevant energy. In addition, it is also connected to a famous conjecture of De Giorgi in early progress, which asserts that entire solutions $u$ to

$$\Delta u = u^3 - u \quad \text{in } \mathbf{R}^n$$

satisfying $|u| \leqslant 1$ and $\partial_{x_n} u > 0$ must be one-dimensional, at least in dimension $n \leqslant 8$. For interested readers, we refer to [GG98] by Ghoussoub–Gui for $n = 2$, [AC00] by Ambrosio–Cabre for $n = 3$, [Sav09] by Savin for $4 \leqslant n \leqslant 8$, and [PKW11] by del Pino–Kowalczyk–Wei for $n \geqslant 9$, for a complete progress of the resolution of the De Giorgi conjecture.

Modica type estimates have now been studied and extended to various situations, for the interested readers, we refer to [CGS94, RR95, FV10, FV11, SS20]. For e.g., it can be applied to the case of systems such as Hénon–Lane–Emden systems. Speaking of this, by writing the fourth-order Hénon equation (1.1) as a system and by a consequence of Phan's pointwise estimate in [Pha12, Lemma 2.7] we know that

$$-\Delta u \geqslant \sqrt{\frac{2}{p+1}} |x|^{\frac{\sigma}{2}} u^{\frac{p+1}{2}} \quad \text{in } \mathbf{R}^n, \tag{1.6}$$

which can be regarded as Modica type estimate in higher order setting. (A simpler version of (1.6) with $\sigma = 0$ already appears in [Sou09], hence (1.6) is sometime referred to Souplet's inequality.) Obviously, to be able to make use of result in [Pha12], the extra condition $0 \leqslant \sigma \leqslant (n-2)(p-1)$ is required. However, it follows from a general Liouville type results in [NY22, NNT24], see also [GG06] and the book [GGS10], that for any $n \geqslant 2$, $p > 1$, and $\sigma \geqslant 0$, it is necessary to have $n \geqslant 5$ and $p \geqslant (n+4+2\sigma)/(n-4)$. Hence, the extra condition holds trivially in the context of the fourth-order Hénon equation (1.1) because

$$0 \leqslant \sigma \leqslant 2\sigma < (n-2)\frac{8+2\sigma}{n-4} \leqslant (n-2)(p-1).$$

In other words, the pointwise estimate (1.6) holds for any $\sigma \geqslant 0$ (with necessarily $n \geqslant 5$ and $p \geqslant (n+4+2\sigma)/(n-4)$).

Concerning to pointwise inequalities for solutions to (1.1), there is another important inequality which has more geometric flavor. Let consider any $C^4$-solution $u$ to (1.1) with $a = 0$ and $p = (n+4)/(n-4)$. Then the conformal metric $g = u^{4/(n-4)} dx^2$ of $\mathbf{R}^n$ has constant Q-curvature $Q_g = 1$. It was expected that the coresponding scalar curvature $R_g$ is also positive (or at least non-negative). By rewriting $g$ as $g = (u^{(n-2)/(n-4)})^{4/(n-2)} dx^2$, we easily compute to get

$$R_g = -\frac{4(n-1)}{n-2} u^{-\frac{n+2}{n-4}} \Delta\big(u^{\frac{n-2}{n-4}}\big) = -\frac{4(n-1)}{n-2} u^{-\frac{n}{n-4}} \Big(\Delta u + \frac{2}{n-4}\frac{|\nabla u|^2}{u}\Big).$$

From this we expect the following pointwise inequality

$$-\Delta u \geqslant \frac{2}{n-4}\frac{|\nabla u|^2}{u} \quad \text{in } \mathbf{R}^n. \tag{1.7}$$

To the best of our knowledge, the pointwise inequality (1.7) should be known among experts. But as far as we know, there is no paper explicity mentioning it, except a very recent paper by Li and Xu [LX25], where they investigate (1.7) for solutions to (1.1) in a more general setting.

Let us now focus on the object of this paper the following interesting pointwise inequality

$$-\Delta u \geqslant \frac{2}{n-4}\frac{|\nabla u|^2}{u} + \sqrt{\frac{2}{p+1-\frac{8}{n(n-4)}}}|x|^{\frac{\sigma}{2}} u^{\frac{p+1}{2}} \quad \text{in } \mathbf{R}^n \tag{1.8}$$

obtained by Fazly et al. in [FWX15] for bounded positive $C^4$-solutions $u$ to (1.1) and for $\sigma \geqslant 0$ but implicitly near 0 and in the supercritical regime $p > (n+4+2\sigma)/(n-4)$. In a sense, (1.8) can be considered as an improvement of Modica type estimate (1.6) and the gradient estimate (1.7).

The argument presented in [FWX15] is nontrivial as it relies on a sophisticated iteration argument in the fashion of the standard Moser proof starting from the rough inequality (1.6). The crown of this delicate argument is a new machinery to obtain pointwise inequalities via iteration, which, as a direct consequence, leads to the existence of positive conformal scalar curvature of the Euclidean space, as mentioned earlier.

Both the estimates (1.6) and (1.7) clearly improve the well-known inequality

$$-\Delta u > 0 \quad \text{in } \mathbf{R}^n, \tag{1.9}$$

which is known as the super-harmonic property for solutions to (1.1). This property was found independent by Lin [Lin98] and Xu [Xu00]; and now plays an important role in the studies of higher-order elliptic equations. See other works [WX99, CAM08, NY22] to name a few.

Although the pointwise inequality (1.8) is remarkable, the statement for its validity apparently has two limitations. The first apparent limitation is the boundedness assumption of solutions, similar to that for the original Modica result. The second limitation is that (1.8) is just known for implicitly small but non-negative $\sigma$. These two limitations come from the approach but with different nature. While the boundedness assumption heavily depends on the use of integral estimates and elliptic estimates, the latter comes from the use of iterations and certain maximum principle argument.

In this paper, we aim to relax these above two limitations. The following theorem is our main result.

**Theorem 1.1** (Fazly–Wei–Xu's inequality). *Let $u$ be a non-negative $C^4$-solution to the fourth-order Hénon equation (1.1) in $\mathbf{R}^n$ with $n \geqslant 5$, $p > 1$, and $\sigma \geqslant 0$. Then, $u$ is strictly positive, $p \geqslant (n+4+2\sigma)/(n-4)$, and the pointwise inequality (1.8) holds, namely*

$$-\Delta u \geqslant \frac{2}{n-4}\frac{|\nabla u|^2}{u} + \sqrt{\frac{2}{p+1-\frac{8}{n(n-4)}}}|x|^{\frac{\sigma}{2}} u^{\frac{p+1}{2}}$$

*holds everywhere in $\mathbf{R}^n$.*

Our theorem above improves the one found by Fazly et al. in [FWX15]. More precisely, no boundedness assumption for solutions to (1.1) is assumed. In addition, the pointwise inequality (1.8) actually holds for any $\sigma \geqslant 0$, not necessary near 0 found in [FWX15]. Moreover, it turns out that the pointwise inequality (1.8) also holds true in the critical case $p = (n+4+2\sigma)/(n-4)$, which was not considered in [FWX15].

We now turn our attention to the necessity of the assumption $\sigma \geqslant 0$ in Theorem 1.1 above. Indeed, on one hand, from the geometric point of view, to be able to define the conformal metric $g = u^{4/(n-4)}dx^2$, it is necessary to have $u \in C^4(\mathbf{R}^n)$, leading to $\sigma \geqslant 0$, thanks to a general regularity result in [GW24]. On the other hand, from the analytic point of view, for $\sigma < 0$, it is clear that our solution space should be $C(\mathbf{R}^n) \cap C^4(\mathbf{R}^n \backslash \{0\})$. Still by the general regularity result in [GW24] we know that if $\sigma$ is sufficiently close to 0, for e.g. $-1 < \sigma < 0$, then $u \in C^2(\mathbf{R}^n)$. Hence, in this scenario the left hand side of (1.8) is defined at 0, but the right hand side is not due to singularity caused by $|x|^{\sigma/2}$.

To prove Theorem 1.1, we use the Bernstein technique in the similar manner of what we did in [NNP18]. Earlier in the [NNP18], we use the Bernstein technique to establish the following pointwise estimate

$$\Delta u \geqslant \frac{1}{2}\frac{|\nabla u|^2}{u} + \sqrt{\frac{2}{p-1-\frac{2}{n}}}\, u^{-\frac{p-1}{2}} \tag{1.10}$$

for positive $C^4$-solutions to $\Delta^2 u = -u^{-p}$ in $\mathbf{R}^n$ with $n \geqslant 3$ and $p > 3$, satisfying some growth condition at infinity. The above pointwise inequality can be considered as a counter-part of (1.8), which also implies the existence of negative conformal scalar curvature of the Euclidean space $\mathbf{R}^n$. In the recent paper [CCR23], still by making use of the Bernstein technique with a new auxiliary function, the authors successfully relax the growth condition at infinity, providing the validity of the weaker version of the above inequality for any positive $C^4$-solutions. It turns out that the new auxiliary function in [CCR23] can be used in the context of (1.1), thus relaxing the boundedness assumption.

Although the two inequalities (1.6) for the power $|x|^{\sigma/2}u^{(p+1)/2}$ and (1.7) for the gradient $u^{-1}|\nabla u|^2$ are well-known, combining them as in (1.8) is not a trivial task even they are scaling invariant. In fact, the gradient term $u^{-1}|\nabla u|^2$ cannot be controlled from the above by the power term $|x|^{\sigma/2}u^{(p+1)/2}$. Let us discuss a simple example. Consider the case $n \geqslant 5$, $\sigma = 0$, and $p = (n+4)/(n-4)$. Then the function

$$U(x) = C(1+|x|^2)^{-(n-4)/2} \tag{1.11}$$

is a solution to (1.1) in $\mathbf{R}^n$ with

$$C = \big((n-4)(n-2)n(n+2)\big)^{\frac{n-4}{8}}. \tag{1.12}$$

Simple calculation shows $(p+3)/4 = (n-2)/(n-4)$ and

$$\nabla U(x) = -(n-4)C(1+|x|^2)^{-(n-2)/2}x,$$

which implies

$$U(x)^{(p+3)/4} = C^{(n-2)/(n-4)}(1+|x|^2)^{-(n-2)/2}$$

and

$$|\nabla U(x)| = (n-4)C(1+|x|^2)^{-(n-2)/2}|x|.$$

Hence,

$$\frac{|\nabla U(x)|}{U(x)^{(p+3)/4}} = (n-4)C^{2/(n-4)}|x|,$$

namely

$$\sup_{\mathbf{R}^n} \frac{|\nabla U|}{U^{(p+3)/4}} = +\infty. \tag{1.13}$$

For general $\sigma > 0$, it is clear that the function $U$ in (1.11) is no longer a solution to (1.1) in $\mathbf{R}^n$. However, it is now standard to note that (1.1) always admits at least one radial solution. In the next result, we confirm (1.13) for any radial solution to (1.1) in $\mathbf{R}^n$ with $\sigma \geqslant 0$.

**Theorem 1.2.** *Let $n \geqslant 5$, $\sigma \geqslant 0$, and $p = (n+4+2\sigma)/(n-4)$. Suppose that $u$ is a radial $C^4$-solution to (1.1) in $\mathbf{R}^n$. Then, there holds*

$$\sup_{\mathbf{R}^n\setminus\{0\}} \frac{|\nabla u|}{|x|^{\sigma/4}u^{(p+3)/4}} = +\infty. \tag{1.14}$$

We prove Theorem 1.2 in section 4. It appears that (1.14) is no longer true in the supercritical regime. In this situation, the quotient $|x|^{-\sigma/4}|\nabla u|u^{-(p+3)/4}$ applied for radial solutions has a finite limit as $|x| \to +\infty$. Details of a proof for this is left for interested readers. Back to (1.14), it is natural to ask whether or not (1.14) remains valid for any $C^4$-solutions, not necessarily radial. We do not know any answer, but it seems that there is an affirmative answer to this question.

Inspired by a question asked in [**FWX15**, Remark 1.4], our next, also last, part of the paper concerns the sharp constants in (1.8), namely the existence of two constants $\alpha_{\mathsf{FWX}}$ and $\beta_{\mathsf{FWX}}$ such that

$$-\Delta u \geqslant \alpha_{\mathsf{FWX}} \frac{|\nabla u|^2}{u} + \beta_{\mathsf{FWX}} |x|^{\frac{\sigma}{2}} u^{\frac{p+1}{2}}$$

holds everywhere in $\mathbf{R}^n$. Thanks to Theorem 1.1, it appears that

$$\alpha_{\mathsf{FWX}} \geqslant \alpha := \frac{2}{n-4} \quad \text{and} \quad \beta_{\mathsf{FWX}} \geqslant \beta := \sqrt{\frac{2}{p+1-\frac{8}{n(n-4)}}}. \tag{1.15}$$

Let us just go back to the radial solution $U$ in the case $\sigma = 0$ given in (1.11). By direct calculation we easily obtain

$$-\Delta U - \frac{2}{n-4}\frac{|\nabla U|^2}{U} = \frac{C(n-4)n}{(1+|x|^2)^{n/2}} > 0$$

and

$$U^{\frac{n}{n-4}} = \frac{C^{n/(n-4)}}{(1+|x|^2)^{n/2}},$$

thanks to $(p+1)/2 = n/(n-4)$. Hence, for any $\varepsilon > 0$ we should have

$$-\Delta U - \Big(\frac{2}{n-4} + \varepsilon\Big)\frac{|\nabla U|^2}{U} = -\frac{C(n-4)(\varepsilon(n-4)|x|^2 - n)}{(1+|x|^2)^{n/2}}.$$

This tells us that the right hand side of the preceding identity becomes negative provided $|x|$ is large enough. This shows that in the case $\sigma = 0$ the coefficient $2/(n-4)$ is sharp, namely $\alpha_{\mathsf{FWX}} = 2/(n-4)$. (In addition, we also have $\beta_{\mathsf{FWX}} = ((n(n-4)/(n^2-4))^{1/2}$, thanks to (1.12).) In fact, we are able to show the sharpness of the coefficient $2/(n-4)$ with arbitrary $\sigma \geqslant 0$. We turn this into a theorem as follows.

**Theorem 1.3** (sharpness). *Let $n \geqslant 5$, $\sigma \geqslant 0$, and $p \geqslant (n+4+2\sigma)/(n-4)$. Then, we have*

$$\alpha_{\mathsf{FWX}} = \frac{2}{n-4}.$$

We prove Theorem 1.3 in section 5 by testing the pointwise inequality (1.8) with radial solutions whose exact behavior at infinity, which is of the form $|x|^{4-n}$, plays an important role. We are successfully in prove the sharpness of the constant $2/(n-4)$ because the two terms $-\Delta u$ and $u^{-1}|\nabla u|^2$ have the same homogeneity in $u$, and the term $|x|^{\sigma/2}u^{(p+1)/2}$ decays faster than the other two terms. On top of these, the crucial ingredient in our proof is that the pointwise inequality (1.8) is invariant under the Kelvin transform, see Lemma 5.1. As far as we know, this feature has not been observed before.

From the above calculation for $U$ corresponding to $\sigma = 0$, it is natural to ask the sharpness of $\beta$ in the case $\sigma > 0$, and we expect that

$$\beta_{\mathsf{FWX}} = \sqrt{\frac{n(n-4)}{(n-2)(n+2)}}.$$

Unfortunately, we do not have an answer. See also Remark 5.5. In the final result, we establish an upper bound for $\beta_{\mathsf{FWX}}$ in the supercritical regime. A precise statement is as follows.

**Theorem 1.4** (upper bound in the supercritical case). *Let $n \geqslant 5$, $\sigma \geqslant 0$, and $p > (n+4+2\sigma)/(n-4)$. Then, there holds*

$$\beta_{\mathsf{FWX}} \leqslant \frac{\theta(n-\theta-2) - \alpha_{\mathsf{FWX}}\theta^2}{\sqrt{\theta(\theta+2)(n-\theta-2)(n-\theta-4)}}$$

*with*

$$\theta = \frac{4+\sigma}{p-1}.$$

We prove Theorem 1.4 in section 6 by testing the pointwise inequality (1.8) with the sharp constants by radial functions. The existence of radial solutions to (1.1) is well-known; see [**NY22**]. Again, the exact behavior at infinity for radial solutions, which is of the form $|x|^{-\theta}$, plays an important role.

Before closing the opening section, we would like to mention that our analysis also works in dimension $n = 3$ leading to some improvement of (1.10) in $\mathbf{R}^3$, compared to known results obtained in [**NNP18**, **CCR23**]. This can also be performed for even higher-order elliptic equations such as the one studied in [**AORW25**, **DN25**]. Details will appear elsewhere.

The organization of this paper is as follows.

## Contents

## 2. Preliminaries

Throughout this section by $u$ we mean a $C^4$-solution to (1.1) in $\mathbf{R}^n$ with $n \geqslant 5$, $p > 1$, and $\sigma \geqslant 0$. At this stage, we just assume $p > 1$. Since our approach is essentially the classical Bernstein technique, the following change of variables

$$v = \log u \quad \text{and} \quad k = \frac{p-1}{2} > 0$$

is very often used. Hence, by tedious calculation one should arrive at the following identities

$$\nabla u = e^v \nabla v = u \nabla v \tag{2.1}$$

and

$$\Delta u = u(|\nabla v|^2 + \Delta v) \tag{2.2}$$

and

$$\Delta^2 u = e^v \big(2|\nabla v|^2 \Delta v + 2\langle \nabla \Delta v, \nabla v\rangle + (\Delta v)^2 + \Delta^2 v + |\nabla v|^4 + 4\nabla^2 v(\nabla v, \nabla v) + \Delta |\nabla v|^2\big), \tag{2.3}$$

where $\nabla^2 v = (v_{ij})$ denotes the Hessian matrix of $v$. Making use of the above three identities we quickly obtain the following lemma.

**Lemma 2.1.** *We have*

$$\begin{aligned}-\Delta^2 v = {} & 2|\nabla v|^2 \Delta v + 4\langle \nabla \Delta v, \nabla v\rangle + (\Delta v)^2 + |\nabla v|^4 \\ & + 4\nabla^2 v(\nabla v, \nabla v) + 2\|\nabla^2 v\|^2 - |x|^\sigma e^{2kv},\end{aligned} \tag{2.4}$$

*where $\|\nabla^2 v\|$ denotes the Hilbert–Schmidt norm on matrices defined to be*

$$\|\nabla^2 v\| = \big(\sum_{i,j} |v_{ij}|^2\big)^{1/2}.$$

*Proof.* Making use of (2.3) and Bochner's identity, namely

$$\Delta |\nabla v|^2 = 2\|\nabla^2 v\|^2 + 2\langle \nabla \Delta v, \nabla v\rangle, \tag{2.5}$$

we arrive at

$$\begin{aligned}u^{-1}\Delta^2 u = {} & 2|\nabla v|^2 \Delta v + 4\langle \nabla \Delta v, \nabla v\rangle + (\Delta v)^2 + \Delta^2 v \\ & + |\nabla v|^4 + 4\nabla^2 v(\nabla v, \nabla v) + 2\|\nabla^2 v\|^2.\end{aligned}$$

Together with $\Delta^2 u = |x|^\sigma u^p$ we obtain

$$\begin{aligned}-\Delta^2 v = {} & 2|\nabla v|^2 \Delta v + 4\langle \nabla \Delta v, \nabla v\rangle + (\Delta v)^2 + |\nabla v|^4 \\ & + 4\nabla^2 v(\nabla v, \nabla v) + 2\|\nabla^2 v\|^2 - |x|^\sigma u^{p-1},\end{aligned}$$

which yields the identity as claimed, thanks to $u^{p-1} = e^{2kv}$. □

Recall the following two universal constants introduced in the preceding section

$$\alpha = \frac{2}{n-4}, \quad \beta = \sqrt{\frac{2}{p+1-\frac{8}{n(n-4)}}},$$

appearing in the desired inequality (1.8). Thanks to (2.1), we know that

$$\alpha\frac{|\nabla u|^2}{u}+\beta|x|^{\frac{\sigma}{2}}u^{\frac{p+1}{2}}=u\big(\alpha|\nabla v|^2+\beta|x|^{\frac{\sigma}{2}}u^{\frac{p-1}{2}}\big)=u(\alpha|\nabla v|^2+\beta|x|^{\frac{\sigma}{2}}e^{kv}).$$

This and (2.2) allow to transform the desired estimate for $u$ into

$$0\geqslant\Delta u+\alpha\frac{|\nabla u|^2}{u}+\beta|x|^{\frac{\sigma}{2}}u^{\frac{p+1}{2}}=u(\Delta v+(\alpha+1)|\nabla v|^2+\beta|x|^{\frac{\sigma}{2}}e^{kv}).$$

This suggests us to define the following auxiliary function $w:\mathbf{R}^n\to\mathbf{R}$ given by

$$w=\Delta v+(\alpha+1)|\nabla v|^2+\beta|x|^{\frac{\sigma}{2}}e^{kv}. \tag{2.6}$$

Our goal is to show that $w\leqslant 0$ everywhere. To this end, we first mention the following weaker version.

**Lemma 2.2.** *We have*

$$w\leqslant\beta|x|^{\frac{\sigma}{2}}e^{kv}\quad\text{in }\mathbf{R}^n. \tag{2.7}$$

*Proof.* From the definition (2.6), it follows that the estimate (2.7) is equivalent to

$$\Delta v+(\alpha+1)|\nabla v|^2\leqslant 0$$

which is further equivalent to

$$u\Delta u+\alpha|\nabla u|^2\leqslant 0\quad\text{in }\mathbf{R}^n,$$

which is nothing but (1.7), which is now well-known. Indeed, as $u\in C^4(\mathbf{R}^n)$, we know that $u$, $\nabla u$, and $\Delta u$ enjoy the following representation

$$u(x)=C(4)\int_{\mathbf{R}^n}\frac{|y|^\sigma u^p(y)}{|x-y|^{n-4}}dy\quad\text{in }\mathbf{R}^n \tag{2.8}$$

and

$$\nabla u(x)=(n-4)C(4)\int_{\mathbf{R}^n}\frac{x-y}{|x-y|^{n-2}}|y|^\sigma u^p(y)dy\quad\text{in }\mathbf{R}^n$$

and

$$-\Delta u(x)=C(2)\int_{\mathbf{R}^n}\frac{|y|^\sigma u^p(y)}{|x-y|^{n-2}}dy\quad\text{in }\mathbf{R}^n$$

with

$$C(4)=\frac{\Gamma(\frac{n-4}{2})}{16\pi^{n/2}},\quad C(2)=\frac{\Gamma(\frac{n-2}{2})}{4\pi^{n/2}}=2(n-4)C(4). \tag{2.9}$$

See [NY22]; see also the proof of Lemma 5.2. In particular, there holds

$$|\nabla u|(x)\leqslant(n-4)C(4)\int_{\mathbf{R}^n}\frac{|y|^\sigma u^p(y)}{|x-y|^{n-3}}dy\quad\text{in }\mathbf{R}^n.$$

Hence, by Hölder's inequality we have

$$\begin{aligned}\alpha|\nabla u|^2&\leqslant 2(n-4)C(4)^2\Big(\int_{\mathbf{R}^n}\frac{|y|^\sigma u^p(y)}{|x-y|^{n-3}}dy\Big)^2\\&\leqslant C(2)C(4)\Big(\int_{\mathbf{R}^n}\frac{|y|^\sigma u^p(y)}{|x-y|^{n-4}}dy\Big)\Big(\int_{\mathbf{R}^n}\frac{|y|^\sigma u^p(y)}{|x-y|^{n-2}}dy\Big)=-u\Delta u.\end{aligned}$$

This completes the proof. □

As mentioned once in Introduction, a more general version of the inequality (1.7) was proved in [LX25] in which the weight $|x|^\sigma$ is replaced by a function $Q$ with appropriate decay at infinity. In the above proof, we offer an alternative argument based on Hölder's inequality. We also note that (2.7) actually holds for any $p > 1$ because the integral representation in [NY22] is available for $p > 1$, but we do not consider subcritical $p$ because no solution exists within this regime.

In turn out that the integral representation (2.8) has other advantages. For e.g. it allows us to estimate solutions from the below. To be more precise, for $|x| \geqslant 1$, we can estimate

$$u(x) \geqslant C(4) \int_{|y| \leqslant 1/2} \frac{|y|^\sigma u^p(y)}{|x-y|^{n-4}} dy \geqslant \Big(\frac{2}{3}\Big)^{n-4} \frac{C(4)}{|x|^{n-4}} \int_{|y| \leqslant 1/2} |y|^\sigma u^p(y) dy,$$

where we have used $|x-y| \leqslant |x| + |y| \leqslant (3/2)|x|$ and $n \geqslant 5$. Hence, there is some $C > 0$ such that

$$u(x) \geqslant \frac{C}{|x|^{n-4}} \quad \text{for } |x| \geqslant 1. \tag{2.10}$$

Remarkably, (2.10) has an important role in our analysis, as we shall use it in the proof of the pointwise inequality (1.8) in dimension $n = 5$.

It turns out that this weak form plays an important role in our analysis. More precisely, instead of starting from Modica type estimate (1.6), which is somewhat weaker than (1.7) as shown in Theorem 1.2, we start with (1.7); see section 3.

Next, we compute and estimate $\Delta w$. This step is very important in Bernstein's method. That said, we have the following claim.

**Lemma 2.3.** *We have*

$$\begin{aligned} \Delta w \geqslant\ & 2(\alpha-1)\langle \nabla v, \nabla w\rangle - I_0 w^2 + I_1 \alpha w |\nabla v|^2 + I_4 \beta |x|^{\frac{\sigma}{2}} w e^{kv} \\ & + I_5 \beta |x|^{\frac{\sigma}{2}} e^{kv} |\nabla v|^2 + \frac{\beta\sigma}{2}\Big(\big(n + \frac{\sigma}{2} - 2\big) + 2(k-\alpha+1)\langle x, \nabla v\rangle\Big)|x|^{\sigma/2-2} e^{kv} \end{aligned} \tag{2.11}$$

*in* $\mathbf{R}^n \backslash \{0\}$, *where* $I_0$, $I_1$, $I_4$, *and* $I_5$ *are constant given in* (2.18). *In addition, if* $\sigma = 0$, *then* (2.11) *holds everywhere in* $\mathbf{R}^n$.

*Proof.* First we consider $\sigma > 0$. Then by working in $\mathbf{R}^n \backslash \{0\}$ we know that any power of $|x|$ is smooth. By the definition of $w$ in (2.6) and Bochner's identity (2.5), we arrive at

$$\begin{aligned} \Delta w &= \Delta^2 v + (\alpha+1)\Delta|\nabla v|^2 + \beta \Delta\big(|x|^{\frac{\sigma}{2}} e^{kv}\big) \\ &= \Delta^2 v + (\alpha+1)\Delta|\nabla v|^2 + \beta \begin{pmatrix} \frac{\sigma}{2}\big(n + \frac{\sigma}{2} - 2\big)|x|^{\frac{\sigma-4}{2}} e^{kv} + k\sigma |x|^{\frac{\sigma-4}{2}} e^{kv} \langle x, \nabla v\rangle \\ + k|x|^{\frac{\sigma}{2}} e^{kv} \Delta v + k^2 |x|^{\frac{\sigma}{2}} e^{kv} |\nabla v|^2 \end{pmatrix} \\ &= \Delta^2 v + 2(\alpha+1)\langle \nabla \Delta v, \nabla v\rangle + 2(\alpha+1)\|\nabla^2 v\|^2 \\ &\quad + \frac{\beta\sigma}{2}\Big(\big(n + \frac{\sigma}{2} - 2\big) + 2k\langle x, \nabla v\rangle\Big)|x|^{\frac{\sigma-4}{2}} e^{kv} + \beta |x|^{\frac{\sigma}{2}}\big(k e^{kv}\Delta v + k^2 e^{kv}|\nabla v|^2\big). \end{aligned}$$

From (2.4) we have

$$\Delta^2 v = -2|\nabla v|^2 \Delta v - 4\langle \nabla \Delta v, \nabla v\rangle - (\Delta v)^2 - |\nabla v|^4 - 4\nabla^2 v(\nabla v, \nabla v) - 2\|\nabla^2 v\|^2 + |x|^\sigma e^{2kv},$$

which then yields

$$\begin{aligned}\Delta w = &-2|\nabla v|^2\Delta v + 2(\alpha-1)\langle\nabla\Delta v,\nabla v\rangle - (\Delta v)^2 - |\nabla v|^4 - 4\nabla^2 v(\nabla v,\nabla v)\\ &+ 2\alpha\|\nabla^2 v\|^2 + |x|^\sigma e^{2kv} + \frac{\beta\sigma}{2}\big((n+\frac{\sigma}{2}-2) + 2k\langle x,\nabla v\rangle\big)|x|^{\frac{\sigma-4}{2}}e^{kv}\\ &+ \beta k^2|x|^{\frac{\sigma}{2}}e^{kv}|\nabla v|^2 + \beta k|x|^{\frac{\sigma}{2}}e^{kv}\Delta v.\end{aligned} \tag{2.12}$$

We now examine the right side of (2.12). First, using (2.6) we rewrite the first term $-|\nabla v|^2\Delta v$ of the right side of (2.12) as

$$-|\nabla v|^2\Delta v = -w|\nabla v|^2 + (\alpha+1)|\nabla v|^4 + \beta|x|^{\frac{\sigma}{2}}e^{kv}|\nabla v|^2. \tag{2.13}$$

The second term $\langle\nabla\Delta v,\nabla v\rangle$ of the right hand side of (2.12) is handled as follows

$$\langle\nabla\Delta v,\nabla v\rangle = \langle\nabla v,\nabla w\rangle - 2(\alpha+1)\nabla^2 v(\nabla v,\nabla v) - k\beta|x|^{\frac{\sigma}{2}}e^{kv}|\nabla v|^2 - \frac{\beta\sigma}{2}|x|^{\frac{\sigma-4}{2}}e^{kv}\langle x,\nabla v\rangle, \tag{2.14}$$

thanks to

$$\nabla\Delta v = \nabla w - (\alpha+1)\nabla|\nabla v|^2 - \beta k|x|^{\frac{\sigma}{2}}e^{kv}\nabla v - \frac{\beta\sigma}{2}|x|^{\frac{\sigma-4}{2}}e^{kv}x$$

and $\langle\nabla v,\nabla|\nabla v|^2\rangle = 2\nabla^2 v(\nabla v,\nabla v)$. The third term $-(\Delta v)^2$ is computed as

$$\begin{aligned}-(\Delta v)^2 = &-w^2 - (\alpha+1)^2|\nabla v|^4 - \beta^2|x|^\sigma e^{2kv} + 2(\alpha+1)w|\nabla v|^2\\ &+ 2\beta|x|^{\frac{\sigma}{2}}we^{kv} - 2(\alpha+1)\beta|x|^{\frac{\sigma}{2}}e^{kv}|\nabla v|^2,\end{aligned} \tag{2.15}$$

and the last term $\beta k|x|^{\frac{\sigma}{2}}e^{kv}\Delta v$ is as follows

$$\beta k|x|^{\frac{\sigma}{2}}e^{kv}\Delta v = \beta k|x|^{\frac{\sigma}{2}}we^{kv} - (\alpha+1)\beta k|x|^{\frac{\sigma}{2}}e^{kv}|\nabla v|^2 - \beta^2 k|x|^\sigma e^{2kv}. \tag{2.16}$$

Combining (2.12)–(2.16) yields

$$\begin{aligned}\Delta w = &-2w|\nabla v|^2 + 2(\alpha+1)|\nabla v|^4 + 2\beta|x|^{\frac{\sigma}{2}}e^{kv}|\nabla v|^2\\ &+\begin{pmatrix}2(\alpha-1)\langle\nabla v,\nabla w\rangle - 4(\alpha^2-1)\nabla^2 v(\nabla v,\nabla v)\\ -2(\alpha-1)k\beta|x|^{\frac{\sigma}{2}}e^{kv}|\nabla v|^2 - (\alpha-1)\beta\sigma|x|^{\frac{\sigma-4}{2}}e^{kv}\langle x,\nabla v\rangle\end{pmatrix}\\ &+\begin{pmatrix}-w^2 - (\alpha+1)^2|\nabla v|^4 - \beta^2|x|^\sigma e^{2kv} + 2(\alpha+1)w|\nabla v|^2\\ +2\beta|x|^{\frac{\sigma}{2}}we^{kv} - 2(\alpha+1)\beta|x|^{\frac{\sigma}{2}}e^{kv}|\nabla v|^2\end{pmatrix}\\ &- |\nabla v|^4 - 4\nabla^2 v(\nabla v,\nabla v) + 2\alpha\|\nabla^2 v\|^2 + |x|^\sigma e^{2kv}\\ &+ \frac{\beta\sigma}{2}\big((n+\frac{\sigma}{2}-2) + 2k\langle x,\nabla v\rangle\big)|x|^{\frac{\sigma-4}{2}}e^{kv} + \beta k^2|x|^{\frac{\sigma}{2}}e^{kv}|\nabla v|^2\\ &+ \beta k|x|^{\frac{\sigma}{2}}we^{kv} - (\alpha+1)\beta k|x|^{\frac{\sigma}{2}}e^{kv}|\nabla v|^2 - \beta^2 k|x|^\sigma e^{2kv}\\ = &\,2\alpha w|\nabla v|^2 - \alpha^2|\nabla v|^4 + 2(\alpha-1)\langle\nabla v,\nabla w\rangle - 4\alpha^2\nabla^2 v(\nabla v,\nabla v) - w^2\\ &+ 2\alpha\|\nabla^2 v\|^2 - \big(\beta^2+\beta^2 k - 1\big)|x|^\sigma e^{2kv} + (k+2)\beta|x|^{\frac{\sigma}{2}}we^{kv}\\ &+ \big(k^2 - (3\alpha-1)k - 2\alpha\big)\beta|x|^{\frac{\sigma}{2}}e^{kv}|\nabla v|^2\\ &+ \frac{\beta\sigma}{2}\big((n+\frac{\sigma}{2}-2) + 2(k-\alpha+1)\langle x,\nabla v\rangle\big)|x|^{\frac{\sigma-4}{2}}e^{kv}.\end{aligned} \tag{2.17}$$

Note that

$$\|\nabla^2 v\|^2 - 2\alpha\nabla^2 v(\nabla v,\nabla v) = \|\nabla^2 v - \alpha\nabla v\otimes\nabla v\|^2 - \alpha^2|\nabla v|^4.$$

Hence,

$$\begin{aligned}\Delta w =\,& 2\alpha w|\nabla v|^2 - \alpha^2(1+2\alpha)|\nabla v|^4 + 2(\alpha-1)\langle \nabla v, \nabla w\rangle - w^2\\ &- \big(\beta^2+\beta^2 k-1\big)|x|^\sigma e^{2kv} + (k+2)\beta|x|^{\frac{\sigma}{2}} w e^{kv} + 2\alpha\|\nabla^2 v - \alpha\nabla v\otimes\nabla v\|^2\\ &+ \big(k^2-(3\alpha-1)k-2\alpha\big)\beta|x|^{\frac{\sigma}{2}} e^{kv}|\nabla v|^2\\ &+ \frac{\beta\sigma}{2}\big((n+\frac{\sigma}{2}-2) + 2(k-\alpha+1)\langle x,\nabla v\rangle\big)|x|^{\frac{\sigma-4}{2}} e^{kv}.\end{aligned}$$

By the Cauchy–Schwarz inequality, we have

$$\begin{aligned}\|\nabla^2 v - \alpha\nabla v\otimes\nabla v\|^2 &\geqslant \frac{1}{n}\big(\Delta v - \alpha|\nabla v|^2\big)^2\\ &= \frac{1}{n}\big(w-(2\alpha+1)|\nabla v|^2 - \beta|x|^{\frac{\sigma}{2}}e^{kv}\big)^2\\ &= \frac{1}{n}\begin{pmatrix} w^2 + (2\alpha+1)^2|\nabla v|^4 + \beta^2|x|^\sigma e^{2kv} - 2(2\alpha+1)w|\nabla v|^2\\ -2\beta|x|^{\frac{\sigma}{2}}we^{kv} + 2(2\alpha+1)\beta|x|^{\frac{\sigma}{2}}|\nabla v|^2 e^{kv}\end{pmatrix}.\end{aligned}$$

Therefore,

$$\begin{aligned}\Delta w \geqslant\,& 2\big(1-\frac{2(2\alpha+1)}{n}\big)\alpha w|\nabla v|^2 - \big(\alpha - \frac{2(2\alpha+1)}{n}\big)(2\alpha+1)\alpha|\nabla v|^4\\ &+ 2(\alpha-1)\langle\nabla v,\nabla w\rangle - \big(1-\frac{2\alpha}{n}\big)w^2 - \big(\beta^2+\beta^2 k - \frac{2\alpha\beta^2}{n} - 1\big)|x|^\sigma e^{2kv}\\ &+ \big(k+2-\frac{4\alpha}{n}\big)\beta|x|^{\frac{\sigma}{2}}we^{kv} + \big(k^2-(3\alpha-1)k-2\alpha+\frac{4\alpha(2\alpha+1)}{n}\big)\beta|x|^{\frac{\sigma}{2}}e^{kv}|\nabla v|^2\\ &+ \frac{\beta\sigma}{2}\big((n+\frac{\sigma}{2}-2)+2(k-\alpha+1)\langle x,\nabla v\rangle\big)|x|^{\frac{\sigma-4}{2}}e^{kv},\end{aligned}$$

namely

$$\begin{aligned}\Delta w \geqslant\,& 2(\alpha-1)\langle\nabla v,\nabla w\rangle - I_0 w^2 + I_1\alpha w|\nabla v|^2 + I_2\alpha|\nabla v|^4 + I_3|x|^\sigma e^{2kv} + I_4\beta|x|^{\frac{\sigma}{2}}we^{kv}\\ &+ I_5\beta|x|^{\frac{\sigma}{2}}e^{kv}|\nabla v|^2 + \frac{\beta\sigma}{2}\big((n+\frac{\sigma}{2}-2)+2(k-\alpha+1)\langle x,\nabla v\rangle\big)|x|^{\frac{\sigma-4}{2}}e^{kv}\end{aligned}$$

with

$$\begin{cases} I_0 = 1-\dfrac{2\alpha}{n} > 0,\\ I_1 = 2\big(1-\dfrac{2(2\alpha+1)}{n}\big) = \dfrac{2(n-6)}{n-4},\\ I_2 = \big(\dfrac{2(2\alpha+1)}{n}-\alpha\big)(2\alpha+1) = 0,\\ I_3 = 1+\dfrac{2\alpha\beta^2}{n} - \beta^2-\beta^2 k = 1-\dfrac{1}{2}\big(p+1-\dfrac{8}{n(n-4)}\big)\beta^2 = 0,\\ I_4 = k+2-\dfrac{4\alpha}{n} = \dfrac{1}{2}\big(p+1-\dfrac{8}{n(n-4)}\big) + I_0,\\ I_5 = k^2-(3\alpha-1)k-2\alpha+\dfrac{4\alpha(2\alpha+1)}{n} = \dfrac{1}{4}\big(p-\dfrac{n+4}{n-4}\big)\big(p+\dfrac{n-8}{n-4}\big) \geqslant 0.\end{cases} \tag{2.18}$$

Hence, we eventually have

$$\begin{aligned}\Delta w \geqslant\,& 2(\alpha-1)\langle\nabla v,\nabla w\rangle - I_0 w^2 + I_1\alpha w|\nabla v|^2 + I_4\beta|x|^{\frac{\sigma}{2}}we^{kv}\\ &+ I_5\beta|x|^{\frac{\sigma}{2}}e^{kv}|\nabla v|^2 + \frac{\beta\sigma}{2}\big((n+\frac{\sigma}{2}-2)+2(k-\alpha+1)\langle x,\nabla v\rangle\big)|x|^{\frac{\sigma-4}{2}}e^{kv}\end{aligned}$$

in $\mathbf{R}^n \setminus \{0\}$, which is the desired estimate as claimed. Finally, the fact that if $\sigma = 0$, then the preceding estimate holds everywhere in $\mathbf{R}^n$ is quite obvious. Therefore, we omit the detail, concluding the present proof. □

For the reader's convenience, we should point out that the condition $p \geqslant (n+4+2\sigma)/(n-4)$ was just used for the very first time to guarantee $I_5 \geqslant 0$.

In the estimate (2.11) above, there is a term involving $\langle x, \nabla v\rangle$ which requires extra care. It is clear that this term appears mainly because of the weight function $|x|^\sigma$. This is also the reason why we have to keep the term involving $I_5$ in the estimate of $\Delta w$ although we know that $I_5 \geqslant 0$.

To go futher, let us introduce a threshold $\sigma^*$ for $\sigma$ defined by

$$\sigma^* = \frac{2(n-2)I_5}{(k-\alpha+1)^2 - I_5}.$$

This will play an important role, as we soon see. But first we show that $\sigma^*$ is well-defined.

**Lemma 2.4.** *We have*

$$(k-\alpha+1)^2 - I_5 > 0.$$

*In particular, there holds* $\sigma^* \geqslant 0$*,*

*Proof.* Note that

$$k-\alpha+1 = \frac{1}{2}\big(p+1-\frac{4}{n-4}\big) = \frac{1}{2}\big(p+\frac{n-8}{n-4}\big).$$

Together with the definition of $I_5$ in (2.18) we have

$$\begin{aligned}(k-\alpha+1)^2 - I_5 &= \frac{1}{4}\big(p+\frac{n-8}{n-4}\big)^2 - \frac{1}{4}\big(p+\frac{n-8}{n-4}\big)\big(p-\frac{n+4}{n-4}\big)\\ &= \frac{1}{2}\big(p+\frac{n-8}{n-4}\big)\frac{n-2}{n-4} > 0\end{aligned}$$

as claimed. □

From the lemma above we deduce that $\sigma^* \geqslant 0$, but we often have $\sigma^* > 0$ except the case when $I_5 = 0$ which occurs only when

$$\sigma = 0 \quad \text{and} \quad p = \frac{n+4}{n-4}.$$

We now describe the role of the threshold $\sigma^*$. In fact, we can further estimate $\Delta w$ from the below provided $0 < \sigma \leqslant \sigma^*$.

**Lemma 2.5.** *Suppose* $0 \leqslant \sigma \leqslant \sigma^*$*. Then, we have*

$$\Delta w \geqslant 2(\alpha-1)\langle \nabla v, \nabla w\rangle - I_0 w^2 + I_1 \alpha w |\nabla v|^2 + I_4 \beta |x|^{\frac{\sigma}{2}} w e^{kv} \tag{2.19}$$

*in* $\mathbf{R}^n \setminus \{0\}$*. In addition, if* $\sigma = 0$*, then* (2.19) *holds everywhere in* $\mathbf{R}^n$*.*

*Proof.* First we observe that if $\sigma = 0$, then the estimate (2.19) is trivial because there is no term involving $\langle x, \nabla v\rangle$ in the estimate (2.11) and $I_5 \geqslant 0$. Hence, we are left with the case $\sigma > 0$. Now under the condition $0 < \sigma \leqslant \sigma^*$ we know that $I_5 > 0$ mainly because $p \geqslant (n+4+2\sigma)/(n-4) > (n+4)/(n-4)$, see (2.18), and that

$$\frac{\beta\sigma}{2}\big(n+\frac{\sigma}{2}-2\big)I_5\beta \geqslant \Big(\frac{\beta\sigma}{2}(k-\alpha+1)\Big)^2$$

because $0 \leqslant \sigma \leqslant \sigma^*$ and $\sigma^*$ is being chosen in such a way that

$$\frac{\beta\sigma^*}{2}(n+\frac{\sigma^*}{2}-2)I_5\beta = \big(\frac{\beta\sigma^*}{2}(k-\alpha+1)\big)^2.$$

Having these facts we can estimate

$$\begin{aligned} I_5\beta|x|^{\frac{\sigma}{2}}e^{kv}|\nabla v|^2 &+ \frac{\beta\sigma}{2}\big((n+\frac{\sigma}{2}-2)+2(k-\alpha+1)\langle x,\nabla v\rangle\big)|x|^{\frac{\sigma-4}{2}}e^{kv} \\ &= I_5\beta|x|^{\frac{\sigma}{2}}e^{kv}\Big|\nabla v + \frac{\beta\sigma}{2}\frac{k-\alpha+1}{I_5\beta}\frac{x}{|x|^2}\Big|^2 \\ &\quad + \Big(\frac{\beta\sigma}{2}(n+\frac{\sigma}{2}-2) - \frac{1}{I_5\beta}\big(\frac{\beta\sigma}{2}(k-\alpha+1)\big)^2\Big)|x|^{\frac{\sigma-4}{2}}e^{kv} \geqslant 0 \end{aligned}$$

in $\mathbf{R}^n\backslash\{0\}$. From this and (2.11) we obtain (2.19) as claimed. □

Having Lemma 2.5 in hand, to obtain the crucial estimate (2.19) for $\Delta w$ from the below, it suffices to examine the inequality $\sigma \leqslant \sigma^*$, Surprisingly, we show that this is always the case.

**Proposition 2.6.** *Under the conditions $n \geqslant 5$, $\sigma > 0$, and $p \geqslant (n+4+2\sigma)/(n-4)$, there always holds $\sigma^* \geqslant 2\sigma > \sigma$. Consequently, the key estimate (2.19) holds for any $\sigma \geqslant 0$ in the sense that it holds in $\mathbf{R}^n\backslash\{0\}$ if $\sigma > 0$ and in $\mathbf{R}^n$ if $\sigma = 0$.*

*Proof.* This follows from direct verification. Indeed, by the definition of $\sigma^*$ we have

$$\sigma^* = \frac{2(n-2)I_5}{(k-\alpha+1)^2 - I_5} = \frac{\frac{1}{2}(n-2)(p-\frac{n+4}{n-4})(p+\frac{n-8}{n-4})}{\frac{1}{2}(p+\frac{n-8}{n-4})\frac{n-2}{n-4}} = \frac{p-\frac{n+4}{n-4}}{\frac{1}{n-4}} \geqslant 2\sigma,$$

where the latter estimate is exactly the the inequality $p \geqslant (n+4+2\sigma)/(n-4)$, which was used for the second time. From this we obtain the key estimate (2.19) for all $\sigma \geqslant 0$. □

*Remark* 2.7. It follows from the proof of Proposition 2.6 above that we only need $p \geqslant (n+4+\sigma)/(n-4)$ for the validity of (2.19) due to $\sigma \leqslant \sigma^*$. This threshold is meaningless due to the known existence result. However, it turns out that this threshold plays some role when working with the Kelvin transform for biLaplacian; see (5.3) in section 5 below.

## 3. Proof of Theorem 1.1

We now prove Theorem 1.1. Let $u$ be a non-negative $C^4$-solution to the fourth-order Hénon equation (1.1) in $\mathbf{R}^n$ with $n \geqslant 5$, $p > 1$, and $\sigma \geqslant 0$. The fact that $u$ is strictly positive and necessarily $p \geqslant (n+4+2\sigma)/(n-4)$ are now well-known; see [**NY22**] and the references therein.

We now prove the pointwise inequality (1.8) by using the Bernstein technique. For clarity, let us start with the simplest case, namely $\sigma = 0$. Since the sign of $I_1$, which depends on the dimension $n$, is important in our argument, we consider two cases $n \geqslant 6$ and $n = 5$ separately.

3.1. **The case $\sigma = 0$ and $n \geqslant 6$.** If $w \leqslant 0$ everywhere in $\mathbf{R}^n$, then there is nothing to prove; see the definition of $w$ in (2.6). Otherwise, there holds

$$\sup_{x\in\mathbf{R}^n} w(x) > 0. \tag{3.1}$$

We shall show that this is not the case. By the continuity of $w$, there is some point $x_0 \in \mathbf{R}^n\backslash\{0\}$ such that $w(x_0) > 0$. Let $\eta : \mathbf{R}^n \to [0,1]$ be a smooth radial cut-off function such that

$$\eta(|x|) = \begin{cases} 0 & \text{if } |x| \geqslant 2, \\ 1 & \text{if } |x| \leqslant 1, \end{cases}$$

and $\eta > 0$ in $B_2$ and that

$$\frac{|\nabla\eta|^2}{\eta} \leqslant C, \quad |\Delta\eta| \leqslant C \quad \text{in } B_2\backslash\overline{B}_1$$

for some constant $C > 0$. Then for arbitrary $R > |x_0| + 1$ but large we let

$$\eta_R(x) = \eta\big(\frac{x}{R}\big)$$

and

$$w_R = \eta_R w.$$

Obviously, there holds

$$\frac{|\nabla\eta_R|^2}{\eta_R} \leqslant \frac{C}{R^2}, \quad |\Delta\eta_R| \leqslant \frac{C}{R^2} \quad \text{in } B_{2R}\backslash\overline{B}_R.$$

As $w_R \equiv 0$ in $\mathbf{R}^n\backslash B_{2R}$, the supremum $\sup_{B_{2R}} w_R$ is achieved at some point $x_R \in B_{2R}$, namely

$$w_R(x_R) = \sup_{B_{2R}} w_R.$$

Thanks to (3.1) and because $w_R \equiv w$ in $B_R$, we may choose $R \gg |x_0| + 1$ in such a way that $w_R(x_R) > 0$. Keep in mind that $\nabla(\eta_R w)(x_R) = 0$. Hence, we have

$$0 \geqslant \eta_R\Delta w + 2\langle\nabla\eta_R, \nabla w\rangle + w\Delta\eta_R = \eta_R\Delta w - 2\frac{|\nabla\eta_R|^2}{\eta_R}w + w\Delta\eta_R \quad \text{at } x_R.$$

Hence

$$\eta_R\Delta w \leqslant \frac{3C}{R^2}w \quad \text{at } x_R, \tag{3.2}$$

which then yields

$$\frac{3C}{R^2}w \geqslant 2(\alpha-1)\langle\nabla v, \nabla w\rangle\eta_R - I_0\eta_R w^2 + I_1\alpha\eta_R w|\nabla v|^2 + I_4\beta\eta_R w e^{kv}$$

at $x_R$, thanks to (2.19). Keep in mind that $w(x_R) > 0$. Unlike other works in the existing literature, in the above estimate, our main term is the linear term of $w$ involving $I_4$, not the quadratic term of $w^2$ involving $I_0$. To make this possible, we use (2.7) to obtain the following crucial estimates

$$\beta e^{kv} \geqslant w > 0 \quad \text{and} \quad \beta w e^{kv} \geqslant w^2 \quad \text{at } x_R.$$

Also, we recall $I_4 = k + 1 - 2\alpha/n + I_0$. From these facts we further have

$$\frac{3C}{R^2}w \geqslant 2(\alpha-1)\eta_R\langle\nabla v, \nabla w\rangle + \Big(k + 1 - \frac{2\alpha}{n}\Big)\eta_R w^2 + 2\big(1 - \frac{2(2\alpha+1)}{n}\big)\alpha\eta_R w|\nabla v|^2$$

at $x_R$, namely

$$\frac{3C}{R^2}w \geqslant -\frac{2(n-6)}{n-4}\eta_R\langle\nabla v, \nabla w\rangle + \frac{4(n-6)}{(n-4)^2}\eta_R w|\nabla v|^2 + \Big(\frac{p+1}{2} - \frac{4}{n(n-4)}\Big)\eta_R w^2 \quad \text{at } x_R. \tag{3.3}$$

Depending on the dimension $n$ we have three cases.

**Case 1**. Suppose $n \geqslant 7$. Keep in mind that $\nabla(\eta_R w)(x_R) = 0$. This and Young's inequality allow us to get

$$\begin{aligned}
-\frac{2(n-6)}{n-4}\langle \nabla v, \nabla w\rangle \eta_R &= \frac{2(n-6)}{n-4}\langle \nabla v, \nabla \eta_R\rangle w \\
&\geqslant -\frac{2(n-6)}{n-4}|\nabla v||\nabla \eta_R| w \\
&\geqslant -\frac{4(n-6)}{(n-4)^2}\eta_R w|\nabla v|^2 - \frac{n-6}{4}\frac{|\nabla \eta_R|^2}{\eta_R} w \quad \text{at } x_R.
\end{aligned}$$

Plugging this into (3.3) gives

$$\frac{(n+6)C}{4R^2} w \geqslant \frac{3C}{R^2} w + \frac{n-6}{4}\frac{|\nabla \eta_R|^2}{\eta_R} w \geqslant \Big(\frac{p+1}{2} - \frac{4}{n(n-4)}\Big)\eta_R w^2 \quad \text{at } x_R.$$

Thus,

$$\Big(\frac{p+1}{2} - \frac{4}{n(n-4)}\Big)^{-1}\frac{(n+6)C}{4R^2} \geqslant (\eta_R w)(x_R) = \sup_{B_{2R}}(\eta_R w) \geqslant w(x_0),$$

thanks to $|x_0| < R$ and $\eta_R \equiv 1$ in $B_R$. Thus, we have just shown that

$$w(x_0) \leqslant \Big(\frac{p+1}{2} - \frac{4}{n(n-4)}\Big)^{-1}\frac{(n+6)C}{4}\frac{1}{R^2}.$$

Letting $R \to +\infty$ we arrive at $w(x_0) \leqslant 0$ which violates $w(x_0) > 0$.

**Case 2**. Suppose $n = 6$. As $p \geqslant 5$, this case is easy to handle because we initially have

$$\frac{3C}{R^2} w \geqslant \frac{8}{3} w^2 \eta_R \quad \text{at } x_R.$$

Thus,

$$\frac{3C}{R^2} \geqslant \frac{8}{3}(\eta_R w)(x_R) = \frac{8}{3}\sup_{B_{2R}}(\eta_R w) \geqslant \frac{8}{3} w(x_0),$$

thanks to $|x_0| < R$ and $\eta_R \equiv 1$ in $B_R$. By arguing similarly as in the case $n \geqslant 7$, we also arrive at contradiction.

3.2. **The case** $\sigma > 0$ **and** $n \geqslant 6$. Next, we consider the case $\sigma > 0$. In this scenario, although our argument is almost the same as in the preceding case, thanks to Proposition 2.6, we still need some extra care due to the presence of the weight $|x|^\sigma$. Indeed, arguing as in the case $\sigma = 0$ we arrive at

$$0 < w_R(x_R) = \sup_{B_{2R}} w_R$$

for some $x_R \in B_{2R}$. This is possible because $w_R \in C(\mathbf{R}^n)$ but not necessary in $C^2(\mathbf{R}^n)$, in fact we can only have $w_R \in C^2(\mathbf{R}^n \backslash \{0\})$, hence we need more job. Keep in mind that $w_R \equiv w$ in $B_R$. Clearly, we have $x_R \neq 0$, otherwise one should have

$$0 < w_R(0) = w(0) = \Delta v(0) + (\alpha+1)|\nabla v|(0)^2 \leqslant 0$$

thanks to (2.7), which is not possible. Thus, by definition, the function $w_R$ is of class $C^2$ in a small neighborhood of $x_R$ in $\mathbf{R}^n \backslash \{0\}$. Then following the remaining argument as in the case $\sigma = 0$ we arrive at

$$\frac{3C}{R^2} w \geqslant 2(\alpha-1)\langle \nabla v, \nabla w\rangle \eta_R - I_0 \eta_R w^2 + I_1 \alpha \eta_R w |\nabla v|^2 + I_4 \beta |x|^{\frac{\sigma}{2}} \eta_R w e^{kv}$$

at $x_R$, see (2.19). Now we make use of

$$\beta |x|^{\frac{\sigma}{2}} e^{kv} \geqslant w > 0 \quad \text{and} \quad \beta |x|^{\frac{\sigma}{2}} w e^{kv} \geqslant w^2 \quad \text{at } x_R$$

to reach (3.3). Hence, the desired inequality $w \leqslant 0$ everywhere follows easily.

3.3. **The case $\sigma = 0$ and $n = 5$.** We again follow the argument used earlier, namely we argue by way of contradiction, hence we obtain $w(x_0) > 0$ for some $x_0 \in \mathbf{R}^n \backslash \{0\}$. Unfortunately, the key estimate (3.3) which now becomes

$$\frac{3C}{R^2} w \geqslant 2\eta_R \langle \nabla v, \nabla w \rangle - 4\eta_R w |\nabla v|^2 + \Big( \frac{p+1}{2} - \frac{4}{5} \Big) \eta_R w^2 \quad \text{at } x_R$$

is not good enough because the term $\eta_R w |\nabla v|^2$ has a wrong sign, and this requires some new idea. Our strategy is to make use of the following change of variable

$$\widehat{w} = e^{\tau v} w \quad \text{in } \mathbf{R}^n$$

for some $\tau \in \mathbf{R}$ to be determined later. (In terms of $u$, the above change of variable is simply $\widehat{w} = u^\tau w$.) Obviously, we have

$$\Delta \widehat{w} = e^{\tau v} \big( w(\tau^2 |\nabla v|^2 + \tau \Delta v) + 2\tau \langle \nabla v, \nabla w \rangle + \Delta w \big).$$

Recall from (2.6) the following $\Delta v = w - (\alpha + 1)|\nabla v|^2 - \beta e^{kv}$. Then, we further have

$$\begin{aligned} e^{-\tau v} \Delta \widehat{w} &= \tau^2 w |\nabla v|^2 + \tau w \big( w - (\alpha+1)|\nabla v|^2 - \beta e^{kv} \big) + 2\tau \langle \nabla v, \nabla w \rangle + \Delta w \\ &= \tau w^2 + \Delta w + 2\tau \langle \nabla v, \nabla w \rangle + \tau(\tau - \alpha - 1) w |\nabla v|^2 - \tau \beta w e^{kv}. \end{aligned}$$

Keep in mind that $I_5 \geqslant 0$. Hence, combining with (2.11) gives

$$\begin{aligned} e^{-\tau v} \Delta \widehat{w} \geqslant{}& 2(\tau + \alpha - 1) \langle \nabla v, \nabla w \rangle - (I_0 - \tau) w^2 \\ &+ \big( I_1 \alpha + \tau(\tau - \alpha - 1) \big) w |\nabla v|^2 + (I_4 - \tau) \beta w e^{kv}. \end{aligned} \tag{3.4}$$

Together with $\alpha = 2$, $I_0 = 1/5$, $I_1 = -2$, and $I_4 = (5p-1)/10$, we thus arrive at

$$e^{-\tau v} \Delta \widehat{w} \geqslant 2(\tau+1) \langle \nabla v, \nabla w \rangle - \big( \frac{1}{5} - \tau \big) w^2 + \big( \tau^2 - 3\tau - 4 \big) w |\nabla v|^2 + \big( \frac{5p-1}{10} - \tau \big) \beta w e^{kv}.$$

Hence, if we choose $\tau = -1$, then we further arrive at

$$e^v \Delta \widehat{w} \geqslant -\frac{6}{5} w^2 + \frac{5p+9}{10} \beta w e^{kv} \geqslant \frac{\beta}{5} w e^{kv} \geqslant \frac{1}{5} w^2 \quad \text{in } \mathbf{R}^n, \tag{3.5}$$

thanks to $p > 1$ and $\beta e^{kv} \geqslant w$. Now we repeat the argument presented in subsection 3.1 applied to $\widehat{w}$ to complete the proof. Indeed, we first observe

$$\sup_{x \in \mathbf{R}^n} \widehat{w}(x) \geqslant e^{-v(x_0)} w(x_0) > 0.$$

Then, we consider

$$\widehat{w}_R = \eta_R \widehat{w},$$

where $\eta_R$ is the cut-off function constructed as before with $R \gg |x_0| + 1$. This leads to the existence of some point $\widehat{x}_R \in B_{2R}$ such that

$$0 < \widehat{w}_R(\widehat{x}_R) = \sup_{B_{2R}} \widehat{w}_R.$$

It is routine to check

$$\frac{3C}{R^2} \widehat{w} \geqslant \eta_R \Delta \widehat{w} \quad \text{at } \widehat{x}_R,$$

see (3.2). Keep in mind that $\widehat{w} = e^{-v} w$. Hence, together with (3.5) we obtain

$$\frac{3C}{R^2} w = \frac{3C}{R^2} e^v \widehat{w} \geqslant e^v \eta_R \Delta \widehat{w} \geqslant \frac{1}{5} \eta_R w^2 \quad \text{at } \widehat{x}_R,$$

namely

$$(\eta_R w)(\widehat{x}_R) \leqslant \frac{15C}{R^2}.$$

Thus, we should have

$$0 < w(x_0) = (\eta_R \widehat{w})(x_0) e^{v(x_0)} \leqslant (\eta_R \widehat{w})(\widehat{x}_R) e^{v(x_0)} = (\eta_R w)(\widehat{x}_R) \frac{e^{v(x_0)}}{e^{v(\widehat{x}_R)}} \leqslant \frac{15C}{R^2} \frac{u(x_0)}{u(\widehat{x}_R)}.$$

Hence, with help of (2.10) and thanks to $\hat{x}_R \in B_{2R}$ we arrive at

$$R^2 u(\hat{x}_R) \geqslant \frac{CR^2}{|\hat{x}_R|} \geqslant \frac{CR^2}{2R} = \frac{CR}{2}$$

for some $C > 0$. Putting the above etimates together we arrive at

$$0 < w(x_0) \leqslant \frac{30Cu(x_0)}{R}$$

for some $C > 0$. This is enough for us to reach contradiction by letting $R \to +\infty$. Thus, $w \leqslant 0$ everywhere. This completes our proof for the case $n = 5$ and $\sigma = 0$.

*Remark* 3.1. The above argument only works in dimension $n = 5$. If $n \geqslant 6$, then it appears that $R^2 u(\hat{x}_R) = O(1)$ as $R \to +\infty$.

3.4. **The case** $\sigma > 0$ **and** $n = 5$. We argue as in the preceding subsection starting from the change of variable

$$\hat{w} = e^{\tau v} w \quad \text{in } \mathbf{R}^n$$

for some $\tau \in \mathbf{R}$ to be specified. Instead of using (2.11) and because $0 < \sigma < \sigma^*$, we use (2.19) to get

$$\begin{aligned} e^{-\tau v} \Delta \hat{w} \geqslant & \, 2(\tau + \alpha - 1)\langle \nabla v, \nabla w \rangle - (I_0 - \tau) w^2 \\ & + \big(I_1 \alpha + \tau(\tau - \alpha - 1)\big) w |\nabla v|^2 c + (I_4 - \tau)\beta w e^{kv}, \end{aligned}$$

which is nothing but the (3.4). Hence, the rest of argument is identical to that presented in the preceding subsection. We omit the details.

## 4. Proof of Theorem 1.2

Let $U$ be any non-negative radial $C^4$-solution to (1.1) in the regime $n \geqslant 5$, $\sigma \geqslant 0$, and $p = (n + 4 + 2\sigma)/(n - 4)$. It is well-known that $U$ is monotone deceasing with respect to the origin. Without loss of generality, with help of simple scaling, we may assume that $U(0) = 1$. We recall that $U$ is a fast decay solution in the sense that

$$C^{-1} r^{4-n} \leqslant U(r) \leqslant C r^{4-n} \quad \text{for all } r \geqslant 1 \tag{4.1}$$

and for some $C \gg 1$; see [**FK19**]. (Keep in mind that in the rest of the paper we shall use (4.1) several times.) Once we have the fast decay property for $U$, we are able to show that (1.14) holds true. Indeed, as

$$C^{-1} r^{4-n} \leqslant U(r) \leqslant C r^{4-n} \quad \text{for all } r \geqslant 1$$

we get

$$C^{-1} r^{-2-\sigma/4} \leqslant r^{\frac{\sigma}{4}} U(r)^{\frac{p-1}{4}} \leqslant C r^{-2-\sigma/4} \quad \text{for all } r \geqslant 1$$

for some new $C \gg 1$. Hence, the integral $\int_0^{+\infty} r^{\sigma/4} U(r)^{(p-1)/4} dr$ converges. Hence

$$\begin{aligned} \sup_{x \in \mathbf{R}^n \setminus \{0\}} \frac{|\nabla U(x)|}{|x|^{\frac{\sigma}{4}} U(x)^{\frac{p+3}{4}}} &= \sup_{r > 0} \frac{-(\log U(r))'}{r^{\frac{\sigma}{4}} U(r)^{\frac{p-1}{4}}} \\ &\geqslant \limsup_{r \to +\infty} \frac{-(\log U(r))'}{r^{\frac{\sigma}{4}} U(r)^{\frac{p-1}{4}}} \\ &\geqslant \limsup_{r \to +\infty} \frac{-\int_0^r (\log U(s))' ds}{\int_0^r s^{\frac{\sigma}{4}} U(s)^{\frac{p-1}{4}} ds} = \limsup_{r \to +\infty} \frac{-\log U(r)}{\int_0^r s^{\frac{\sigma}{4}} U(s)^{\frac{p-1}{4}} ds}, \end{aligned}$$

where we have used $U(0) = 1$ to obtain the last estimate. Keep in mind that

$$0 < \int_0^r s^{\frac{\sigma}{4}} U(s)^{\frac{p-1}{4}} ds < \int_0^{+\infty} s^{\frac{\sigma}{4}} U(s)^{\frac{p-1}{4}} ds$$

for all $r > 0$ and

$$-\log U(r) = \log \frac{1}{U(r)} \geqslant \log(C^{-1} r^{n-4}) = -\log C + (n-4)\log r \to +\infty$$

as $r \to +\infty$. This tells us

$$\sup_{x \in \mathbf{R}^n \setminus \{0\}} \frac{|\nabla U(x)|}{|x|^{\frac{\sigma}{4}} U(x)^{\frac{p+3}{4}}} = +\infty$$

as claimed. Our proof for (1.14) in the regime $p = (n+4+2\sigma)/(n-4)$ is complete.

*Remark* 4.1. It is natural to ask the value of

$$\inf_{x \in \mathbf{R}^n \setminus \{0\}} \frac{|\nabla U(x)|}{|x|^{\frac{\sigma}{4}} U(x)^{\frac{p+3}{4}}}.$$

In the case $\sigma = 0$, the above infimum must be zero because $U \in C^4$ and $U(0) > 0$. This is still true if $0 < \sigma \leqslant 12$. For e.g. if $0 < \sigma \leqslant 4$, then we have

$$0 = \lim_{r \searrow 0} r^{1-\frac{\sigma}{4}} \Delta U(r) = \lim_{r \searrow 0} \frac{(r^{n-1}U'(r))'}{r^{n-2+\frac{\sigma}{4}}} = \big(n-1+\frac{\sigma}{4}\big) \lim_{r \searrow 0} \frac{U'(r)}{r^{\frac{\sigma}{4}}}$$

by the l'Hopital rule. Here, we have used $U \in C^2(\mathbf{R}^n)$. For $4 < \sigma \leqslant 12$, the above argument is not enough, but as $U \in C^4(\mathbf{R}^n)$, we further have

$$\begin{aligned} 0 = \lim_{r \searrow 0} r^{3-\frac{\sigma}{4}} \Delta^2 U(r) &= \big(n-3+\frac{\sigma}{4}\big) \lim_{r \searrow 0} \frac{(\Delta U)'(r)}{r^{-2+\frac{\sigma}{4}}} \\ &= \big(n-2+\frac{\sigma}{4}\big)\big(-1+\frac{\sigma}{4}\big) \lim_{r \searrow 0} \frac{\Delta U(r)}{r^{-1+\frac{\sigma}{4}}} \end{aligned}$$

still by the l'Hopital rule. For general $\sigma > 12$, the answer seems to be complicated and heavily depends on the regularity of $U$ near the origin.

## 5. Proof of Theorem 1.3

In this section, we prove the sharpness of the constant $2/(n-4)$ in the pointwise inequality (1.8) in the supercritical case $p \geqslant (n+4+2\sigma)/(n-4)$. This is done by testing (1.8) with radial solutions to (1.1).

First, we let $U$ be any positive radial $C^4$-solution to (1.1) in the regime $n \geqslant 5$, $\sigma \geqslant 0$, and $p \geqslant (n+4+2\sigma)/(n-4)$. For convenience, it is no harm to assume that $U$ is radially symmetric with respect to the origin. In this scenario, we already know that $U$ is strictly decreasing and there is some $c > 0$ such that

$$U(r) \leqslant c r^{-\frac{4+\sigma}{p-1}}$$

and

$$|U^{(k)}(r)| < c r^{-\frac{4+\sigma}{p-1} - k} \quad \text{for } 1 \leqslant k \leqslant 4$$

in $(0, +\infty)$; see [NN26]. By Theorem 1.1 we already know that

$$-\Delta U \geqslant \frac{2}{n-4} \frac{|\nabla U|^2}{U} + \beta |x|^{\frac{\sigma}{2}} U^{\frac{p+1}{2}} \quad \text{in } \mathbf{R}^n.$$

To prove $\alpha_{\mathsf{FWX}} = 2/(n-4)$ we argue by way of contradiction that

$$\alpha_{\mathsf{FWX}} = \frac{2}{n-4} + \varepsilon$$

for some $\varepsilon > 0$, namely

$$\Delta U + \big(\frac{2}{n-4} + \varepsilon\big) \frac{|\nabla U|^2}{U} + \beta |x|^{\frac{\sigma}{2}} U^{\frac{p+1}{2}} \leqslant 0 \quad \text{in } \mathbf{R}^n. \tag{5.1}$$

Consider the following Kelvin transform for biLaplacian with center at the origin

$$\widehat{U}(x) = \frac{1}{|x|^{n-4}} U\Big(\frac{x}{|x|^2}\Big) \quad \text{in } \mathbf{R}^n \backslash \{0\}.$$

Then $\widehat{U}$ solves

$$\Delta^2 \widehat{U} = |x|^{\widehat{\sigma}} \widehat{U}^p \quad \text{in } \mathbf{R}^n \backslash \{0\} \tag{5.2}$$

with $\widehat{\sigma} = (n-4)p - (n+4+\sigma)$. In addition, $\widehat{U}$ is radially symmetric with respect to the origin. Then, under the supercritical condition and $\sigma \geqslant 0$, there holds

$$\frac{n+\widehat{\sigma}}{n-4} < p < \frac{n+4+2\widehat{\sigma}}{n-4} \quad \text{and} \quad \widehat{\sigma} \geqslant 0. \tag{5.3}$$

To go further, we need the following lemma.

**Lemma 5.1.** *There holds*

$$\Delta \widehat{U}(x) + \frac{2}{n-4} \frac{|\nabla \widehat{U}|^2}{\widehat{U}}(x) = \frac{1}{|x|^n} \Big[ \Delta U\Big(\frac{x}{|x|^2}\Big) + \frac{2}{n-4} \frac{|\nabla U|^2}{U}\Big(\frac{x}{|x|^2}\Big) \Big] \tag{5.4}$$

*and*

$$|x|^{\frac{\widehat{\sigma}}{2}} \widehat{U}^{\frac{p+1}{2}}(x) = \frac{1}{|x|^n} \Big|\frac{x}{|x|^2}\Big|^{\frac{\sigma}{2}} U^{\frac{p+1}{2}}\Big(\frac{x}{|x|^2}\Big) \tag{5.5}$$

*everywhere in* $\mathbf{R}^n \backslash \{0\}$.

*Proof.* This can be checked directly. Indeed, (5.5) is simple because

$$|x|^{\frac{\widehat{\sigma}}{2}} \widehat{U}^{\frac{p+1}{2}}(x) = |x|^{\frac{1}{2}[(n-4)p-(n+4+\sigma)]} \Big(|x|^{4-n} U\Big(\frac{x}{|x|^2}\Big)\Big)^{\frac{p+1}{2}} = \frac{1}{|x|^n} \Big|\frac{x}{|x|^2}\Big|^{\frac{\sigma}{2}} U^{\frac{p+1}{2}}\Big(\frac{x}{|x|^2}\Big).$$

We now verify (5.4). For convenience, we denote $W = U^{\frac{n-2}{n-4}}$. For convenience, we denote

$$\breve{W}(x) = \frac{1}{|x|^{n-2}} W\Big(\frac{x}{|x|^2}\Big)$$

(also known as the Kelvin transform for Laplacian). Then, there holds $\breve{W} = \widehat{U}^{\frac{n-2}{n-4}}$. Clearly,

$$(\Delta \widehat{U}^{\frac{n-2}{n-4}})(x) = \Delta \breve{W}(x) = \frac{1}{|x|^{n+2}} (\Delta W)\Big(\frac{x}{|x|^2}\Big) = \frac{1}{|x|^{n+2}} (\Delta U^{\frac{n-2}{n-4}})\Big(\frac{x}{|x|^2}\Big).$$

Then,

$$\widehat{U}^{\frac{2}{n-4}}(x)\Big(\Delta \widehat{U} + \frac{2}{n-4} \frac{|\nabla \widehat{U}|^2}{\widehat{U}}\Big)(x) = \frac{1}{|x|^{n+2}} U^{\frac{2}{n-4}}\Big(\frac{x}{|x|^2}\Big)\Big(\Delta U + \frac{2}{n-4} \frac{|\nabla U|^2}{U}\Big)\Big(\frac{x}{|x|^2}\Big).$$

From this we obtain (5.4) because

$$\widehat{U}^{\frac{2}{n-4}}(x) = \frac{1}{|x|^2} U^{\frac{2}{n-4}}\Big(\frac{x}{|x|^2}\Big).$$

This completes the proof. □

Using (5.4)–(5.5) we can transfer (5.1) into

$$\Delta \widehat{U}(x) + \frac{2}{n-4} \frac{|\nabla \widehat{U}|^2}{\widehat{U}}(x) + \beta |x|^{\frac{\widehat{\sigma}}{2}} \widehat{U}^{\frac{p+1}{2}}(x) + \frac{\varepsilon}{|x|^n} \frac{|\nabla U|^2}{U}\Big(\frac{x}{|x|^2}\Big) \leqslant 0 \tag{5.6}$$

in $\mathbf{R}^n \backslash \{0\}$. Keep in mind that $U$ is radial. Notice that

$$\frac{|\nabla \widehat{U}|^2}{\widehat{U}}(x) = \frac{1}{|x|^n} \frac{|\nabla U|^2}{U}\Big(\frac{x}{|x|^2}\Big) + \frac{2(n-4)}{|x|^{n-1}} U'\Big(\frac{1}{|x|}\Big) + \frac{(n-4)^2}{|x|^{n-2}} U\Big(\frac{x}{|x|^2}\Big).$$

Hence, combining with (5.6) we further arrive at

$$\begin{aligned}\Delta \widehat{U}(x)+&\Big(\frac{2}{n-4}+\varepsilon\Big)\frac{|\nabla \widehat{U}|^2}{\widehat{U}}(x)+\beta|x|^{\frac{\widehat{\sigma}}{2}}\widehat{U}^{\frac{p+1}{2}}(x)\\ &-\frac{2(n-4)\varepsilon}{|x|^{n-1}}U'\Big(\frac{1}{|x|}\Big)-\frac{(n-4)^2\varepsilon}{|x|^{n-2}}U\Big(\frac{x}{|x|^2}\Big)\leqslant 0\end{aligned} \tag{5.7}$$

in $\mathbf{R}^n\backslash\{0\}$. Keep in mind that

$$\widehat{U}(x)\sim\frac{U(0)}{|x|^{n-4}}\quad\text{near }+\infty$$

thanks to $U\in C^4(\mathbf{R}^n)$. For convenience, we denote $\tau=n-4$ and use the Emden–Fowler transformation

$$V(t)=e^{\tau t}\widehat{U}(e^t)\quad\text{in }\mathbf{R}.$$

Then, it is known that $V$ has a positive limit $U(0)$ as $t\to+\infty$, say $L>0$. Then as

$$|\nabla\widehat{U}(x)|=-\widehat{U}'(r)=-\big(V'(t)-\tau V(t)\big)r^{-\tau-1}$$

and

$$\Delta\widehat{U}(x)=\widehat{U}''(r)+\frac{n-1}{r}\widehat{U}'(r)=\big(V''(t)+(n-2\tau-2)V'(t)-\tau(n-\tau-2)V(t)\big)r^{-\tau-2}$$

it follows from (5.7) that

$$\begin{aligned}0\geqslant &V(t)V''(t)+(n-2\tau-2)V(t)V'(t)-\tau(n-\tau-2)V(t)^2\\ &+\Big(\frac{2}{n-4}+\varepsilon\Big)\big(V'(t)-\tau V(t)\big)^2+\beta e^{-\tau\frac{p+1}{2}t}V(t)^{\frac{p+3}{2}}\\ &-2(n-4)\varepsilon e^{-(n-1)t}U'(e^{-t})-(n-4)^2\varepsilon e^{-(n-2)t}U(e^{-t})\end{aligned} \tag{5.8}$$

after canceling out the term $r^{-\tau-3}$ from both sides. Now we assume for a mommet that

$$V^{(k)}(t)\to 0\quad\text{as }t\to+\infty \tag{5.9}$$

for $k=1,2$. Then by passing the above inequality to the limit as $t\to+\infty$ we arrive at

$$\tau(n-\tau-2)L^2\geqslant\Big(\frac{2}{n-4}+\varepsilon\Big)(\tau L)^2.$$

which gives

$$\frac{2}{n-4}+\varepsilon\leqslant\frac{n-\tau-2}{\tau}=\frac{2}{n-4}.$$

This is not true. Thus, we must have $\alpha_{\mathsf{FWX}}=2/(n-4)$ as claimed.

In the rest of the present proof, we verify (5.9), namely we show that $V^{(k)}(t)\to 0$ as $t\to+\infty$ with $k=1,2$. Since the argument is dedicated, we state and prove as lemmas. Before doing so, we start with some useful information for $\widehat{U}$.

**Lemma 5.2.** *We have the following*

$$\widehat{U}'<0,\quad\Delta\widehat{U}<0,\quad(\Delta\widehat{U})'>0$$

*in* $(0,+\infty)$. *In addition, the function*

$$r\mapsto r\widehat{U}'(r)+(n-4)\widehat{U}(r)$$

*is non-negative and non-increasing in* $(0,+\infty)$.

*Proof.* First, by making use of [**NY22**, Proposition 1.3] applied to the equation (5.2) together with (5.3) we know that $\widehat{U}$ solves

$$\widehat{U}(x) = C(4)\int_{\mathbf{R}^n} \frac{|y|^{\widehat{\sigma}}\widehat{U}^p(y)}{|x-y|^{n-4}}dy \quad \text{in } \mathbf{R}^n\backslash\{0\}. \tag{5.10}$$

The lower bound for $p$ in (5.3) is to ensure that $\widehat{U}$ is a distributional solution to (5.2). Then (5.10) follows because $\widehat{U}$ enjoys the so-called ring condition. Reasoning similarly and with help of [**NY22**, Lemma 2.4] we also know that $-\Delta\widehat{U}$ is a distributional solution to the second-order equation $-\Delta w = |x|^{\widehat{\sigma}}\widehat{U}^p$. Hence, it solves

$$-\Delta\widehat{U}(x) = C(2)\int_{\mathbf{R}^n} \frac{|y|^{\widehat{\sigma}}\widehat{U}^p(y)}{|x-y|^{n-2}}dy \quad \text{in } \mathbf{R}^n\backslash\{0\}. \tag{5.11}$$

Here, the coefficients $C(2)$ and $C(4)$ are already given in (2.9). In particular, there hold

$$-\Delta\widehat{U} > 0$$

in $(0,+\infty)$. Now as $\widehat{U}$ and $-\Delta\widehat{U}$ are non-negative and superharmonic in $\mathbf{R}^n\backslash\{0\}$, the remaining two inequalities follow from standard arguments; see e.g. [**NNT26**, Remark A.2(3)].

Next, we verify the monotonicity of the function $r\widehat{U}'(r)+(n-4)\widehat{U}(r)$ by using an idea from [**CMM93**]. First, we observe

$$(r(-\Delta\widehat{U})'(r)+(n-2)(-\Delta\widehat{U})(r))' = r\big((-\Delta\widehat{U})''(r)+\frac{n-1}{r}(-\Delta\widehat{U})'(r)\big) = -r\Delta^2\widehat{U} < 0$$

in $(0,+\infty)$. Hence $r(-\Delta\widehat{U})'(r)+(n-2)(-\Delta\widehat{U})(r)$ is strictly decreasing. If there were some $r_0>0$ such that $r_0(-\Delta\widehat{U})'(r_0)+(n-2)(-\Delta\widehat{U})(r_0)<0$, then we would have

$$r(-\Delta\widehat{U})'(r)+(n-2)(-\Delta\widehat{U})(r) < r_0(-\Delta\widehat{U})'(r_0)+(n-2)(-\Delta\widehat{U})(r_0) =: C < 0$$

in $[r_0,+\infty)$, namely

$$(-\Delta\widehat{U})'(r) < (-\Delta\widehat{U})'(r)+\frac{n-2}{r}(-\Delta\widehat{U})(r) < \frac{C}{r} < 0.$$

Integrating over $[r_0,r]$ gives

$$-\infty < -(-\Delta\widehat{U})(r_0) < (-\Delta\widehat{U})(r)-(-\Delta\widehat{U})(r_0) < C\log\frac{r}{r_0}$$

in $[r_0,+\infty)$. Sending $r\nearrow+\infty$ we obtain contradiction, hence

$$r(-\Delta\widehat{U})'+(n-2)(-\Delta\widehat{U}) \geqslant 0 \tag{5.12}$$

everywhere in $(0,+\infty)$. From this we are able to obtain

$$\begin{aligned}
-(r^{n-1}\widehat{U}')'' &= \big(r^{n-1}(-\Delta\widehat{U})\big)' \\
&= (n-1)r^{n-2}(-\Delta\widehat{U})+r^{n-1}(-\Delta\widehat{U})' \\
&= r^{n-2}\big(r(-\Delta\widehat{U})'+(n-2)(-\Delta\widehat{U})\big)+r^{n-2}(-\Delta\widehat{U}) \\
&\geqslant r^{n-2}(-\Delta\widehat{U}).
\end{aligned}$$

Thus,

$$-(r^{n-1}\widehat{U}')'' \geqslant r^{-1}r^{n-1}(-\Delta\widehat{U}) = r^{-1}\big(-r^{n-1}\widehat{U}'\big)'$$

which futher implies

$$\big(r(r^{n-1}\widehat{U}')'-2r^{n-1}\widehat{U}'\big)' \leqslant 0 \quad \text{in } (0,+\infty).$$

Hence, $r(r^{n-1}\widehat{U}')' - 2r^{n-1}\widehat{U}'$ is non-increasing in $(0,+\infty)$. Notice that

$$\begin{aligned} r^{n-1}|\widehat{U}'(r)| &= r^{n-1}\big|\big(r^{4-n}U(r^{-1})\big)'\big| \\ &\leqslant (n-4)r^2U(r^{-1}) + r|U'(r^{-1})| \\ &\leqslant cr^{2+\frac{4+\sigma}{p-1}} \end{aligned}$$

and

$$\begin{aligned} r|(r^{n-1}\widehat{U}'(r))'| &= r\big|\big(r^{n-1}\big(r^{4-n}U(r^{-1})\big)'\big)'\big| \\ &\leqslant (n-4)r^2U(r^{-1}) + (n-5)rU'(r^{-1}) + U''(r^{-1}) \\ &\leqslant cr^{2+\frac{4+\sigma}{p-1}} \end{aligned}$$

for some $c>0$. Putting the above estimates together we arrive at

$$\lim_{r\searrow 0}\big(r(r^{n-1}\widehat{U}')' - 2r^{n-1}\widehat{U}'\big) = 0$$

yielding

$$r(r^{n-1}\widehat{U}')' \leqslant 2r^{n-1}\widehat{U}' \quad \text{in } (0,+\infty),$$

or equivalently

$$(n-1)r^{n-1}\widehat{U}' + r^n\widehat{U}'' \leqslant 2r^{n-1}\widehat{U}' \quad \text{in } (0,+\infty).$$

Hence,

$$(n-3)\widehat{U}' + r\widehat{U}'' \leqslant 0 \quad \text{in } (0,+\infty).$$

Thus

$$(r\widehat{U}'(r) + (n-4)\widehat{U})' \leqslant 0 \quad \text{in } (0,+\infty),$$

proving the monotonicity of $r\widehat{U}' + (n-4)\widehat{U}$ in $(0,+\infty)$ as claimed. Finally, the inequality

$$\widehat{U}'(r) + (n-4)\widehat{U}(r) \geqslant 0 \quad \text{in } (0,+\infty) \tag{5.13}$$

can be checked by a similar way of proving (5.12). □

Let us now start with the limit of $V'$.

**Lemma 5.3.** *There holds $V'(t)\to 0$ as $t\to+\infty$.*

*Proof.* Since $V(t)$ has a limit as $t\to+\infty$, the conclusion follows by making use of Proposition A.1 provided $V''$ is bounded in $[1,+\infty)$. A proof of the boundedness of $V''$ in $[1,+\infty)$ is actually provided in the proof of the next lemma. So in the rest of the present proof, we offer an alternative argument. First we recall the identity

$$r\widehat{U}'(r) = \big(V'(t) - \tau V(t)\big)r^{-\tau}.$$

Together with $\widehat{U}(r) = r^{-\tau}V(t)$ we arrive at

$$r^\tau\big(r\widehat{U}'(r) + (n-4)\widehat{U}(r)\big) = V'(t) + (n-4-\tau)V(t) = V'(t).$$

Hence, to conclude the limit of $V'$, it suffices to show

$$\lim_{r\to+\infty} r^{n-4}\big(r\widehat{U}'(r) + (n-4)\widehat{U}(r)\big) = 0. \tag{5.14}$$

Thanks to Lemma 5.2, we have

$$\big(r\widehat{U}'(r) + (n-4)\widehat{U}(r)\big)' \leqslant 0 \quad \text{for all } r>0.$$

Resolving the preceding inequality gives

$$r\big(\widehat{U}''(r) + \frac{n-3}{r}\widehat{U}'(r)\big) \leqslant 0 \quad \text{for all } r>0,$$

namely

$$(r^{n-3}\widehat{U}'(r))' = r^{n-3}\big(\widehat{U}''(r) + \frac{n-3}{r}\widehat{U}'(r)\big) \leqslant 0 \quad \text{for all } r>0.$$

Thus, we have just shown that the function $r \mapsto r^{n-3}\widehat{U}'(r)$ is monotone decreasing in $(0,+\infty)$. As $U \in C^4$ and

$$0 \leqslant r^{3-n}|\widehat{U}'(r)| \leqslant (n-4)U(r^{-1}) + r^{-1}|U'(r^{-1})| \tag{5.15}$$

in $(0,+\infty)$ we conclude that the function $r \mapsto r^{n-3}\widehat{U}'(r)$ is bounded in $[1,+\infty)$. This together with its monotonicity concludes that the limit

$$\lim_{r\to+\infty} r^{n-3}\widehat{U}'(r)$$

exists. From this we can apply the l'Hôpital rule to get

$$(n-4)\lim_{r\to+\infty} r^{n-4}\widehat{U}(r) = -\lim_{r\to+\infty} r^{n-3}\widehat{U}'(r),$$

namely there holds

$$\lim_{r\to+\infty}\big(r^{n-3}\widehat{U}'(r) + (n-4)r^{n-4}\widehat{U}(r)\big) = 0.$$

From this we immediately get (5.14) because

$$0 \leqslant V'(t) = r^{n-4}\big(r\widehat{U}'(r) + (n-4)\widehat{U}(r)\big) = r^{n-3}\widehat{U}'(r) + (n-4)r^{n-4}\widehat{U}(r).$$

The proof is complete. □

Now we study the limit of $V''$. We shall prove the following.

**Lemma 5.4.** *There holds $V''(t) \to 0$ as $t \to +\infty$.*

*Proof.* By Lemma 5.3 we know that $V'(t) \to 0$ as $t \to +\infty$. Thus, we have the conclusion by making use of Proposition A.1 provided $V^{(3)}$ is bounded in $[1,+\infty)$. In the rest of the proof we prove that $V^{(3)}$ is indeed bounded in $[1,+\infty)$. Clearly, $V$ is bounded in $[1,+\infty)$ and it follows from Lemma 5.3 that $V'$ is also bounded in $[1,+\infty)$. In fact, the boundedness of $V'$ can also be realized from the identity

$$V'(t) = \tau V(t) + r^{\tau+1}\widehat{U}'(r) \tag{5.16}$$

and the boundedness of $r^{\tau+1}\widehat{U}'(r)$, which was already proved in (5.15). By differentiating (5.16) we arrive at

$$V''(t) = (2\tau+1)V'(t) - \tau(\tau+1)V(t) + r^{\tau+2}\widehat{U}''(r), \tag{5.17}$$

for which we can conclude the boundedness of $V''$ provided $r^{\tau+2}\widehat{U}''(r)$ is bounded. To verify the boundedness of $r^{\tau+2}\widehat{U}''(r)$ we first use

$$0 \geqslant \Delta\widehat{U}(r) = \widehat{U}''(r) + (n-1)r^{-1}\widehat{U}'(r) \quad \text{for } r>0$$

to obtain the estimate

$$r^{\tau+2}|\widehat{U}''(r)| \leqslant (n-1)r^{\tau+1}|\widehat{U}'(r)| - r^{\tau+2}\Delta\widehat{U}(r) \quad \text{for } r>0.$$

Hence, it suffices to bound $-r^{\tau+2}\Delta\widehat{U}(r)$ from the above, which can be argued as follows. First, we observe from $(r^{n-1}\widehat{U}'(r))' = r^{n-1}\Delta\widehat{U}$, $\widehat{U}'<0$, and $(\Delta\widehat{U})' \geqslant 0$, the

following

$$r^{n-1}\widehat{U}'(r) \leqslant r^{n-1}\widehat{U}'(r) - \widehat{U}'(1) = \int_1^r s^{n-1}\Delta\widehat{U}(s)ds$$
$$\leqslant \Delta\widehat{U}(r)\int_1^r s^{n-1}ds = \frac{r^n}{n}\Delta\widehat{U}(r) \leqslant 0$$

for $r \geqslant 1$. Keep in mind the definition of $\widehat{U}$ and (5.13), hence, we quickly have

$$Cr^{-\tau} \geqslant \widehat{U}(r) \geqslant -\frac{1}{n-2}r\widehat{U}'(r) \geqslant -\frac{1}{n(n-2)}r^2\Delta\widehat{U}(r) \geqslant 0$$

for $r \geqslant 1$. Thus, we have just shown that the boundedness of $-r^{\tau+2}\Delta\widehat{U}(r)$ on $[1,+\infty)$ as claimed. Finally, by differentiating the identity (5.17) we obtain

$$V^{(3)}(t) = 3(\tau+1)V''(t) - (3\tau^2+6\tau+2)V'(t) + (\tau^3+3\tau^2+2\tau)V(t) + r^{\tau+3}\widehat{U}^{(3)}(r).$$

Therefore, to bound $V^{(3)}$, it suffices to bound $r^{\tau+3}\widehat{U}^{(3)}(r)$, and this can be seen by making use of

$$(\Delta\widehat{U})'(r) = \widehat{U}^{(3)}(r) + (n-1)r^{-1}\widehat{U}''(r) - (n-1)r^{-2}\widehat{U}'(r) \quad \text{for } r > 0,$$

which yields

$$0 \leqslant r^{\tau+3}|\widehat{U}^{(3)}(r)| \leqslant (n-1)r^{\tau+2}|\widehat{U}''(r)| + (n-1)r^{\tau+1}|\widehat{U}'(r)| + r^{\tau+3}(\Delta\widehat{U})'(r).$$

Then, we make use of (5.12) to get

$$0 \leqslant r^{\tau+3}(\Delta\widehat{U})'(r) \leqslant -(n-2)r^{\tau+2}\Delta\widehat{U}(r) \quad \text{for } r > 0.$$

From this we obtain the boundedness of $r^{\tau+3}(\Delta\widehat{U})'(r)$ in $[1,+\infty)$, so is $V^{(3)}$. □

We have a remark before closing this section.

*Remark* 5.5. It is natural to ask the sharpness of $\beta$, and we expect, for e.g. in the critical case $p = (n+4+2\sigma)/(n-4)$, that

$$\beta_{\mathsf{FWX}} = \sqrt{\frac{n(n-4)}{(n-2)(n+2)}}.$$

To look for $\beta_{\mathsf{FWX}}$, one needs to go back to (5.8) to get

$$-V(t)V''(t) + (n-2)V(t)V'(t) \geqslant \frac{2}{n-4}(V'(t))^2 + \beta e^{-(2+\frac{\sigma}{2})t}V(t)^{\frac{p+3}{2}}$$

thanks to $\tau(n-\tau-2)V^2 = \alpha_{\mathsf{FWX}}(\tau V)^2$. As both sides of the above inequality decay to zero, to look for the sharp constant $\beta_{\mathsf{FWX}}$ one needs to know the exact further decay rate of $V$, $V'$, and $V''$ at infinity.

## 6. Proof of Theorem 1.4

In Theorem 1.3, we obtain the sharpness of the constant $\alpha_{\mathsf{FWX}}$, which is equal to $2/(n-4)$, in the pointwise inequality (1.8) in the case $p \geqslant (n+4+2\sigma)/(n-4)$. Unfortunately, the proof does not work for the case of the constant $\beta_{\mathsf{FWX}}$. In this section, we provide an upper bound for the sharp constant $\beta_{\mathsf{FWX}}$ in the inequality (1.8). To illustrate our idea, we also provide an upper bound for the sharp constant $\alpha_{\mathsf{FWX}}$.

Let $U$ be any positive radial $C^4$-solution to (1.1) in the regime $n \geqslant 5$, $\sigma \geqslant 0$, and $p > (n+4+2\sigma)/(n-4)$. In this scenario, $U$ is slow decay solution, see [NN26]. By Theorem 1.1 we know that

$$-\Delta U \geqslant \alpha \frac{|\nabla U|^2}{U} + \beta |x|^{\frac{\sigma}{2}} U^{\frac{p+1}{2}}$$

in $\mathbf{R}^n$. Here the constants $\alpha$ and $\beta$ are given in (1.15). Via the Emden–Fowler transformation, we let $|x| = r$ and define

$$V(t) = e^{\theta t} U(e^t) \quad \text{in } \mathbf{R}$$

namely $|x| = r = e^t$ or $t = \log r$, with $\theta = (4+\sigma)/(p-1)$. Equivalently, we have

$$U(r) = r^{-\theta} V(\log r) \quad \text{in } (0, +\infty).$$

Keep in mind that $V > 0$ everywhere. Then it is well-known that

$$|\nabla U(x)| = U'(r) = \big(V'(t) - \theta V(t)\big) r^{-(\theta+1)}$$

and that

$$\Delta U(x) = U''(r) + \frac{n-1}{r} U'(r) = \big(V''(t) + (n-2\theta-2)V'(t) - \theta(n-\theta-2)V(t)\big) r^{-(\theta+2)}.$$

Hence, resolving the inequality

$$-U\Delta U \geqslant \alpha_{\mathsf{FWX}} |\nabla U|^2 + \beta_{\mathsf{FWX}} |x|^{\frac{\sigma}{2}} U^{\frac{p+3}{2}}$$

gives

$$\begin{aligned} -V(t)V''(t) - (n-2\theta-2)V(t)V'(t) + \theta(n-\theta-2)V(t)^2 \\ - \alpha_{\mathsf{FWX}} \big(V'(t) - \theta V(t)\big)^2 \geqslant \beta_{\mathsf{FWX}} V(t)^{\frac{p+3}{2}}. \end{aligned}$$

Sending $t \to +\infty$ yields

$$\frac{\theta(n-\theta-2) - \alpha_{\mathsf{FWX}} \theta^2}{\sqrt{\theta(\theta+2)(n-\theta-2)(n-\theta-4)}} \geqslant \beta_{\mathsf{FWX}}.$$

This completes the proof in the regime $p > (n+4+2\sigma)/(n-4)$.

## Acknowledgments

The first author would like to thank Professor Xingwang Xu and Doctor Mingxiang Li for useful discussion on this topic, especially during their visit to TSIMF in Sanya in July 2025.

## Appendix A. A calculus lemma

In this appendix we mention a simple and perhaps well-known calculus lemma that we have used in the proof of Theorem 1.3.

**Proposition A.1.** *Let $f : [1, +\infty) \to \mathbf{R}$ be a $C^2$-function such that $f''$ is bounded and the limit $\lim_{t\to+\infty} f(t)$ exists. Then, we have*

$$\lim_{t\to+\infty} f'(t) = 0.$$

*Proof.* Subtracting from $f$ by a suitable constant we can assume from the beginning that $f(t) \to 0$ as $t \to +\infty$. Let $M > 0$ be such that $|f''(t)| \leqslant M$ for all $t \geqslant 1$. Making use of the Taylor expansion gives

$$f(t+h) = f(t) + f'(t)h + \frac{1}{2}f''(t+\lambda h)h^2$$

for some $\lambda \in (0,1)$ depending on $t$ and $h$. From this we quickly get

$$|f'(t)| \leqslant \frac{|f(t+h) - f(t)|}{h} + \frac{Mh}{2}.$$

In view of the limit $f(t) \to 0$ as $t \to +\infty$, given arbitrary $\varepsilon \in (0,1)$ we can find some $t_0 \gg 1$ in such a way that

$$|f(t)| \leqslant \frac{\varepsilon}{2} \quad \text{for all} t \geqslant t_0.$$

Putting the above estimates together yields

$$|f'(t)| \leqslant \frac{\varepsilon}{h} + \frac{Mh}{2} \quad \text{for all} t \geqslant t_0.$$

Now we choose $h = \sqrt{2\varepsilon/M}$ to obtain

$$0 \leqslant |f'(t)| \leqslant \sqrt{2\varepsilon M} \quad \text{for all} t \geqslant t_0.$$

From this we obtain the desired limit. □

(Q.A. Ngô) University of Science, Vietnam National University, Hanoi, Vietnam, ORCID iD: 0000-0002-3550-9689

*Email address*: nqanh@vnu.edu.vn

(N.T. Trung) University of Science, Vietnam National University, Hanoi, Vietnam

*Email address*: trungisp58@gmail.com